\documentclass[11pt]{article}%
\usepackage{amsmath,amssymb,color}
\usepackage{amsmath}
\usepackage{amsfonts}
\usepackage{amssymb}
\usepackage{graphicx}%
\providecommand{\U}[1]{\protect\rule{.1in}{.1in}}
\definecolor{mycolorred}{rgb}{1, 0, 0}

\def\<{\langle}
\def\>{\rangle}

\newtheorem{theorem}{Theorem}[section]

\newtheorem{corollary}[theorem]{Corollary}

\newtheorem{lemma}[theorem]{Lemma}

\newtheorem{proposition}[theorem]{Proposition}
\newtheorem{remark}[theorem]{Remark}

\numberwithin{equation}{section}
\begin{document}

\title{Upper bounds for the function solution of the homogeneous $2D$ Boltzmann
equation with hard potential}
\author{\textsc{Vlad Bally}\thanks{Universit\'e Paris-Est, LAMA (UMR CNRS, UPEMLV,
UPEC), INRIA, F-77454 Marne-la-Vall\'ee, France. Email:
\texttt{bally@univ-mlv.fr} }\smallskip}
\date{}
\maketitle

\begin{abstract}
We deal with $f_{t}(dv),$ the solution of the homogeneous $2D$ Boltzmann
equation without cutoff. The initial condition $f_{0}(dv)$ may be any
probability distribution (except a Dirac mass). However, for sufficiently hard
potentials, the semigroup has a regularization property (see \cite{[BF]}):
$f_{t}(dv)=f_{t}(v)dv$ for every $t>0.$ The aim of this paper is to give upper
bounds for $f_{t}(v),$ the most significant one being of type $f_{t}(v)\leq
Ct^{-\eta}e^{-\left\vert v\right\vert ^{\lambda}}$ for some $\eta,\lambda>0.$

\end{abstract}


\parindent 0pt

\medskip


\medskip


\bigskip

\section{Introduction and main results}

We are concerned with the solution of the two dimensional homogenous Boltzmann
equation:
\begin{equation}
\partial_{t}f_{t}(v)=\int_{\mathbb{R}^{2}}dv_{\ast}\int_{-\pi/2}^{\pi
/2}d\theta\left\vert v-v_{\ast}\right\vert ^{\gamma}b(\theta)(f_{t}(v^{\prime
})f_{t}(v_{\ast}^{\prime})-f_{t}(v)f_{t}(v_{\ast})).\label{bo1}%
\end{equation}
Here $f_{t}(v)$ is a non-negative measure on $\mathbb{R}^{2}$ which represents
the density of particles having velocity $v$ in a model for a gas in dimension
two, and, with $R_{\theta}$ being the rotation of angle $\theta,$
\[
v^{\prime}=\frac{v+v_{\ast}}{2}+R_{\theta}\left(  \frac{v-v_{\ast}}{2}\right)
,\quad v_{\ast}^{\prime}=\frac{v+v_{\ast}}{2}-R_{\theta}\left(  \frac
{v-v_{\ast}}{2}\right)  .
\]
The function $b:[-\frac{\pi}{2},\frac{\pi}{2}]\diagdown\{0\}\rightarrow
\mathbb{R}$ will be assumed to satisfy the following hypothesis:%
\begin{align}
(H_{\nu})\quad i)\quad\exists0 &  <c<C\quad s.t.\quad c\left\vert
\theta\right\vert ^{-(1+\nu)}\leq b(\theta)\leq C\left\vert \theta\right\vert
^{-(1+\nu)}\label{bo1a}\\
\quad ii)\quad\forall k &  \in N,\exists C_{k}\quad s.t.\quad\left\vert
b^{(k)}(\theta)\right\vert \leq C_{k}\left\vert \theta\right\vert
^{-(k+1+\nu)}.\nonumber
\end{align}
In \cite{FM} it is proved that, for every $\nu\in(0,\frac{1}{2})$ and
$\gamma\in(0,1],$ the above equation has a unique weak solution. More
precisely: under the assumption $(H_{\nu})$ and the integrability condition
$\int e^{\left\vert v\right\vert ^{\lambda}}f_{0}(dv)<\infty$ for some
$\lambda\in(\gamma,2),$\cite{FM} shows that there exists a unique weak
solution $f_{t}$ of (\ref{bo1}) which starts from $f_{0}.$ Furthermore the
solution satisfies $\sup_{t\leq T}\int e^{\left\vert v\right\vert
^{\lambda^{\prime}}}f_{t}(dv)<\infty$ for every $\lambda^{\prime}<\lambda.$
Throughout the paper these hypotheses are in force.

Below, $f_{0}(dv)$ will be a probability distribution which is not assumed to
be absolutely continuous with respect to the Lebesgue measure (we will just
assume that $f_{0}(dv)$ is not a Dirac mass $\delta_{v_{0}}(dv)$ - in this
trivial case the corresponding solution is $f_{t}(dv)=f_{0}(dv)=\delta_{v_{0}%
}(dv)$ for every $t>0).$

Our first aim is to give sufficient conditions under which, for every $t>0,$
$f_{t}(dv)$ is absolutely continuous and to study the regularity of its
density $f_{t}(v)$: see Theorem \ref{L} below. This problem has already been
addressed in \cite{[BF]} for the same equation and under the same assumptions,
in \cite{F1} for the three dimensional Boltzmann equation and in \cite{ADVW}
for the Boltzmann equation in arbitrary dimension (however, in this last
paper, $f_{0}(dv)$ is assumed to be absolutely continuous and to have finite
entropy). In the case of Maxwell molecules (where $\gamma=0)$ this problem is
addressed in \cite{GM} and \cite{F3}, and for Landau equation in \cite{GMN}
.\ These last three papers are the pioneering papers concerning the
probabilistic approach to the regularity problem.

Our second aim is to give upper bounds for $f_{t}(v):$ see Corollary \ref{R}
and Theorem \ref{Bound}\ bellow. This is actually the main contribution of the
paper. Two aspects of our bounds are noteworthy:

\begin{itemize}
\item for any $t>0,v\mapsto f_{t}(v)$ decays exponentially fast (spatial
exponential decay) and

\item $t\mapsto f_{t}(v)$ blows up at most polynomially as $t\rightarrow0$
(blow-up in time).
\end{itemize}

Under a cutoff condition, and if the initial value is a function which is
upper bounded by a Maxwellian potential, bounds on the spatial decay of
$f_{t}$ were proved in \cite{GPV}. Our result applies even when the initial
condition is a measure (not necessarily absolutely continuos) and we work with
the equation without cutoff. In \cite{K} the author discusses upper bounds of
polynomial type for an initial condition which is a smooth function. This is
rather different from our framework, as our bounds are exponential and no
regularity of the initial condition is required.

In order to precisely state our results, some notation is required. We denote
by $\left\Vert \cdot\right\Vert _{p}$ respectively by $\left\Vert
\cdot\right\Vert _{q,p}$ the norm in $L^{p}$ respectively in the Sobolev space
$W^{q,p}$ on $\mathbb{R}^{2}$. For $p>1$ we denote by $p_{\ast}$ the conjugate
of $p.$ We fix $\nu\in(0,\frac{1}{2})$ and $\gamma\in(0,1],$ we suppose that
$f_{0}(dv)$ is not a Dirac mass, and that for some $\lambda\in(\gamma,2)$ one
has $\int e^{\left\vert v\right\vert ^{\lambda}}f_{0}(dv)<\infty.$

We consider a non decreasing function $\rho:\mathbb{R}_{+}\rightarrow
\mathbb{R}_{+}$ such that $\rho(u)=1$ for $u\in(0,1),\rho(u)=u$ for
$u\in(2,\infty)$ and $\rho\in C^{\infty}(\mathbb{R}_{+})$ and, for
$0<\lambda^{\prime}<\lambda,$ we define
\begin{equation}
\Phi_{\lambda^{\prime}}(v)=e^{\rho(\left\vert v\right\vert ^{\lambda^{\prime}%
})}. \label{Lp2'}%
\end{equation}
The important point is that $\Phi_{\lambda^{\prime}}(v)=e^{\left\vert
v\right\vert ^{\lambda^{\prime}}}$ for $\left\vert v\right\vert \geq2;$ the
function $\rho$ is used just to avoid singularities of the derivatives of
$\Phi_{\lambda^{\prime}}$ around $v=0.$ The specific choice of $\rho$ impacts
just the constants $C$ (which anyway are not explicit).

Moreover we will use an auxiliary function $\varphi:[0,\infty)\rightarrow
\lbrack0,\infty)$ defined by
\begin{equation}
\varphi(\alpha)=\frac{(1-\nu)(1+\gamma+\alpha)}{1+\nu(\gamma+\alpha
)}-1\label{Lp1}%
\end{equation}
and we denote by $\alpha_{\ast}$ the unique solution of the equation
$\varphi(\alpha_{\ast})=\alpha_{\ast}$ (see (\ref{alpha}) for the explicit
value of $\alpha_{\ast}).$ We also denote%
\begin{equation}
\eta=\frac{2(\varphi(2)-1)}{\varphi(2)-2}\left(  \frac{13(1+\alpha_{\ast
})(2+\nu)}{\nu}-1\right)  .\label{Lp3b}%
\end{equation}
We need to impose $\eta>0.$\ Direct computation shows that%
\begin{align}
\varphi(0) &  >0\quad\Leftrightarrow\quad\nu<\frac{\gamma}{2\gamma
+1}\label{Lp3e}\\
\varphi(2) &  >2\quad\Leftrightarrow\quad\nu<\frac{\gamma}{4\gamma
+9}\label{Lp3d}%
\end{align}

\begin{theorem}
\label{L}\textbf{A}. Let us assume that $\nu<\frac{\gamma}{2\gamma+1}.$ Then
the measure $f_{t}(dv)$ is absolutely continuous with respect to the Lebesgue
measure on $\mathbb{R}^{2}$. We denote by $f_{t}(v)$ its density (that is
$f_{t}(dv)=f_{t}(v)dv).$

Let $\lambda^{\prime}<\lambda.$

\textbf{B \ a}. If $\nu<\frac{\gamma}{4\gamma+9}$ then $\Phi_{\lambda^{\prime
}}f_{t}\in L^{p}(\mathbb{R}^{2})$ for every $p>1$ and (with $\eta$ given in
(\ref{Lp3b}))%
\begin{equation}
\left\Vert \Phi_{\lambda^{\prime}}f_{t}\right\Vert _{p}\leq\frac{C}{t^{\eta}}.
\label{Lp3c}%
\end{equation}
\textbf{b}. If $\frac{\gamma}{4\gamma+9}\leq\nu<\frac{\gamma}{2\gamma+1}$ then
$\alpha_{\ast}<2$ and $\Phi_{\lambda^{\prime}}f_{t}\in L^{p}(\mathbb{R}^{2})$
for every $1<p<\frac{2}{2-\alpha_{\ast}}.$

\textbf{C}.\textbf{a} If $\nu<\frac{\gamma}{4\gamma+9}$ then $\Phi
_{\lambda^{\prime}}f_{t}\in W^{q,p}(\mathbb{R}^{2})$ for $q=1,2$ and
$1<p<p_{q},$ with%
\begin{equation}
p_{1}=\frac{2(1+\nu(\gamma+2))}{1-\gamma+11\nu+5\nu\gamma}\quad and\quad
p_{2}=\frac{2(1+\nu(\gamma+2))}{2-\gamma+13\nu+6\nu\gamma}. \label{Lp4}%
\end{equation}
Moreover for every $p<p_{q}$ one has
\begin{equation}
\left\Vert \Phi_{\lambda^{\prime}}f_{t}\right\Vert _{q,p}\leq\frac{C}{t^{\eta
}}. \label{Lp5}%
\end{equation}

\textbf{b}. If $\frac{\gamma}{4\gamma+9}\leq\nu<\frac{\gamma}{3\gamma+4}$ then
$\Phi_{\lambda^{\prime}}f_{t}\in W^{1,p}(\mathbb{R}^{2})$ for every
$1<p<\frac{2}{3-\alpha_{\ast}}.$
\end{theorem}

We stress that the precise power $\eta$ in $t^{-\eta}$ in (\ref{Lp3c}) and
(\ref{Lp5}) is due to the technical approach that we use. We do not expect it
to be optimal (see the point $D$ in Lemma \ref{A} for more precise estimates,
which themselves are not optimal). However this guarantees that the blow up of
$f_{t}$ as $t\rightarrow0$ is at most polynomial.

In order to be able to compare this result with the ones in the papers which
we quoted before, take $s>1$ and $\nu=\frac{2}{s-1},\gamma=\frac{s-5}{s-1}:$
these are the values which are significant in the case of the three
dimensional Boltzmann equation. Our condition $\gamma>0$ implies that $s>5$;
in the literature this case is known as the "hard potential" case. With\ this
choice of $\nu$ and of $\gamma$ we have $\nu<\frac{\gamma}{2\gamma+1}$ iff
$s>9$ and $\nu<\frac{\gamma}{4\gamma+9}$ iff $s>16+\sqrt{193}\sim30.$ The
regularity results of the above theorem are analogous with the ones in
\cite{[BF]}, though not identical. In \cite{F1} one deals with the real three
dimensional equation (without cutoff) and obtains absolute continuity for a
larger range for $s$ then in the above theorem. However the $L^{p}$ estimates
obtained in our paper are stronger: we obtain $\Phi_{\lambda^{\prime}}f_{t}\in
L^{p}(\mathbb{R}^{2})$ instead of $f_{t}\in L^{2}(\mathbb{R}^{2}).$ Moreover,
we obtain bounds depending polynomially on $t\downarrow0$. The result of
\cite{ADVW} is stronger in the sense that it applies to equations in any
dimension, but it requires that the initial condition is already a function
(so it is not really possible to compare them).

We give now some consequences of the previous result concerning the tails of
$f_{t}(dv):$

\begin{corollary}
\label{R}Suppose that $\nu<\frac{\gamma}{2\gamma+1}.$ For every $\lambda
^{\prime}<\lambda$ there exists a constant $C\geq1$ (depending on
$\lambda^{\prime})$ such that for every $R>1,t\in(0,1]$%
\begin{align}
f_{t}(\{v  &  :\left\vert v\right\vert \geq R\})\leq\frac{C}{t^{\kappa}%
}e^{-R^{\lambda^{\prime}}}\quad with\label{Lp6}\\
\kappa &  =\frac{13(2+\nu)(1-\nu)(1+\gamma)}{\nu(1+\nu\gamma)}-1. \label{Lp6a}%
\end{align}

\end{corollary}

We give now the upper bound for $f_{t}(v):$

\begin{theorem}
\label{Bound}Suppose that $\nu<\frac{\gamma}{4\gamma+9}.$ Then $p_{1}>2$
(given in (\ref{Lp4}))\ and $f_{t}\in C^{0,\chi}$ (H\"{o}lder continuous
functions of order $\chi)$ with $\chi=1-\frac{2}{p_{1}}$ for all $t>0$.
Moreover for every $\lambda^{\prime}<\lambda$ there exists $C\geq1$ such that%
\begin{equation}
\left\vert f_{t}(v)\right\vert \leq\frac{C}{t^{\eta}}e^{-\left\vert
v\right\vert ^{\lambda^{\prime}}} \label{bo2a}%
\end{equation}
with $\eta$ given in (\ref{Lp3b}). Finally, there exists $C\geq1$ such that
for every $v,w\in\mathbb{R}^{2}$ with $\left\vert w-v\right\vert \leq1$%
\begin{equation}
\left\vert f_{t}(w)-f_{t}(v)\right\vert \leq\frac{C}{t^{\eta}}e^{-\left\vert
v\right\vert ^{\lambda^{\prime}}}\left\vert w-v\right\vert ^{\chi}.
\label{bo2bc}%
\end{equation}

\end{theorem}

To our best knowledge, both the time-space estimate (\ref{bo2a}) and the
H\"{o}lder continuity of  $f_{t}$ and the estimate (\ref{bo2bc}) are new. For
the Landau equation, lower and upper bounds for the solution have been
obtained in \cite{GMN} using integration by parts techniques based on the
classical Malliavin calculus.\ This approach is not directly possible in our
framework because of the singularities that appear in the problem.\ We will
use similar but ultimately different techniques below.

Corollary \ref{R}\ and Theorem \ref{Bound} are the main contributions of our
paper. The drawback of our approach is that it applies only to "very hard
potentials" ($s>9$ for (\ref{Lp6}) respectively $s>30$ for (\ref{bo2a}) and
(\ref{bo2bc})). Moreover, the exponent $\eta$ in the polynomial blow-up
$t^{-\eta}$ is not expected to be optimal.

The proofs are based on a "balance argument" which is interesting in itself,
and may be useful in other settings. We summarize it below.

Consider a family of random variables $F_{\varepsilon}\sim f_{\varepsilon
}(v)dv,\varepsilon>0$ and a random variable $F.$ Suppose that $F_{\varepsilon
}-F\rightarrow0$ and $f_{\varepsilon}\uparrow\infty$ as $\varepsilon
\rightarrow0$, in a certain sense. If the convergence to zero is sufficiently
faster then the blow-up\ of $f_{\varepsilon}$, then one is able to prove that
the law of $F$ \ has a density and to obtain some regularity of the density.
This idea first appears in \cite{[FP]} and has been used ever since in several
papers (see \cite{[De]} for example). In these papers the "balance" between
the speed of convergence to zero and the blow-up is built by using Fourier
analysis. Later on, in \cite{DR} the authors introduced a new method based on
a Besov space criterion, which turns out to be significantly more powerful
then the one based on Fourier analysis. This is the method used in \cite{F1}
in the case of the three dimensional Boltzmann equation (see also \cite{DF}).
Finally, in \cite{BCpreprint}, a third method which is close to interpolation
theory was introduced.

The criterion that we use in the present paper is an improvement of the latter
method: we give an abstract framework in which an integration by parts formula
can be applied and we quantify the blow-up of $f_{\varepsilon}$ in terms of
the weights appearing in the corresponding integration by parts formula. Let
us be more precise.

Consider a family of random variables $F_{\varepsilon}$ with values in
$\mathbb{R}^{d}$ and $G_{\varepsilon}$ with values in $[0,1],\varepsilon>0$.
Associate to them the measures $\mu_{\varepsilon}$ given by%
\[
\int\varphi d\mu_{\varepsilon}=E(\varphi(F_{\varepsilon})G_{\varepsilon}).
\]
The random variables $G_{\varepsilon}$ play a technical role, and will be used
in some localization procedure.

We assume that for every $\varepsilon>0$ and every multi-index $\alpha$ one
may find a random variable $H_{\alpha,\varepsilon}$ such that the following
integration by parts formula holds:%
\begin{equation}
E(\partial^{\alpha}\varphi(F_{\varepsilon})G_{\varepsilon})=E(\varphi
(F_{\varepsilon})H_{\alpha,\varepsilon})\quad\forall\varphi\in C_{b}^{\infty
}(\mathbb{R}^{d}). \label{ip1}%
\end{equation}
Here $\alpha=(\alpha_{1},...,\alpha_{m})\in\{1,...,d\}^{m},$ is a multi-index
of length $\left\vert \alpha\right\vert =m$ and $\partial^{\alpha}$ is the
derivative associated to $\alpha.$

Additionally, we assume that $H_{\alpha,\varepsilon}$ may be chosen such that,
for every $q\in\mathbb{N}$ and $p>1$, there exist some constants $\widehat
{H}_{q,p}$ and $a,b,\varepsilon_{\ast}\geq0$ such that for every
$0<\varepsilon<\varepsilon_{\ast}$
\begin{equation}
\sup_{\left\vert \alpha\right\vert \leq q}\left\Vert H_{\alpha,\varepsilon
}\right\Vert _{p}\leq\widehat{H}_{q,p}\varepsilon^{-b(q+a)}. \label{ip2}%
\end{equation}
In particular this implies that $\mu_{\varepsilon}(dv)=f_{\varepsilon}(v)dv$
with $f_{\varepsilon}\in C^{\infty}(\mathbb{R}^{d}).$

Moreover, we consider a random variable $F\in R^{d}$ and we assume that there
exists $\beta>0$ and $C_{\ast}\geq1$ such that%
\begin{equation}
\left\Vert G_{\varepsilon}-1\right\Vert _{2}+\left\Vert F_{\varepsilon
}-F\right\Vert _{1}\leq C_{\ast}\varepsilon^{\beta}. \label{ip3}%
\end{equation}

Finally we consider a function $\Phi:\mathbb{R}^{d}\rightarrow\mathbb{R}_{+}$
which belongs to $C^{\infty}(\mathbb{R}^{d}),$ is convex and there exists
$C\geq1$ such that $\Phi(x)\geq\frac{1}{C}\Phi(y)$ if $\left\vert x\right\vert
\geq\left\vert y\right\vert .$ Moreover we assume that for each $h\in
\mathbb{N}$ and for each multi-index $\alpha$ there exist some constants
$c_{1},c_{2}$ (depending on $h$ and $\alpha)$ such that
\begin{equation}
(1+\left\vert x\right\vert )^{h}(1+\left\vert \partial^{\alpha}\Phi
(x)\right\vert )\leq c_{1}\Phi^{c_{2}}(x).\label{ip4}%
\end{equation}
The typical examples are $\Phi(x)=(1+\left\vert x\right\vert ^{2})^{r}$ and
$\Phi_{\lambda}$ defined in (\ref{Lp2'}). We denote by $\mathcal{C}$\ the
class of these functions and for $\Phi\in\mathcal{C}$ and $\theta\geq0$ we
denote%
\begin{equation}
C_{\theta}(\Phi)=E(\Phi^{\theta}(F))+\sup_{\varepsilon>0}E(\Phi^{\theta
}(F_{\varepsilon})).\label{ip5}%
\end{equation}
Our criteria are the following:

\begin{theorem}
\label{Criterion}\textbf{A}. Let $F\in R^{d}$ be a random variable. Suppose
that one is able to find a family $F_{\varepsilon}\in\mathbb{R}^{d}$ and
$G_{\varepsilon}\in\mathbb{R},\varepsilon>0$ which verify (\ref{ip1}),
(\ref{ip2}) and (\ref{ip3}). Fix $q\in\mathbb{N}$ and $p>1$ and assume that
(recall that $p_{\ast}$ is the conjugate of $p)$%
\begin{equation}
\beta>b(1+q+\frac{d}{p_{\ast}}).\label{ip6}%
\end{equation}
Then $P(F\in dx)=f(x)dx$ with $f\in W^{q,p}(\mathbb{R}^{d}).$

\textbf{B}. Consider a function $\Phi\in\mathcal{C}$ such that $C_{\theta
}(\Phi)<\infty$ for every $\theta>0.$ Let $q\in\mathbb{N}$ and $p>1,\delta>0$
be given. Assume that (\ref{ip6}) holds. There exist some constants $C\geq
1$,$\theta\geq1$ and $h_{\ast}\geq1$ (depending on $q,d,\beta,b,p$ and
$\delta)$ such that for $h\geq h_{\ast}$ one has%
\begin{align}
\left\Vert \Phi f\right\Vert _{q,p}  &  \leq\Gamma_{\Phi,\theta}(q,h,p)\quad
with\label{ip7}\\
\Gamma_{\Phi,\theta}(q,h,p)  &  :=C\times(C_{\ast}+C_{\theta}(\Phi
))\times\left(  h^{2b}\widehat{H}_{2h+q+d,p_{\ast}}^{1/2h}\right)
^{(1+\delta)(1+q+d/p_{\ast})} \label{ip8}%
\end{align}
with $C_{\ast}$ given in (\ref{ip3}), $\widehat{H}_{2h+q+d,p_{\ast}}$\ given
in\ (\ref{ip2}) and $C_{\theta}(\Phi)$ given in (\ref{ip5}).

\textbf{C}. Suppose that (\ref{ip6}) holds for $q=1$ and $p>d.$ Then $f\in
C^{0,\chi}(\mathbb{R}^{d})$ with $\chi=1-\frac{d}{p}$\ and we have
\begin{equation}
\left\vert f(x)\right\vert \leq\frac{1}{\Phi(x)}\times\left\Vert \Phi
f\right\Vert _{1,p}\leq\frac{1}{\Phi(x)}\times\Gamma_{\Phi,\theta}(1,h,p).
\label{ip9}%
\end{equation}
Moreover, let%
\begin{equation}
\widehat{\nabla}\Phi(x)=\sup_{\left\vert x-y\right\vert \leq1}\left\vert
\nabla\Phi(y)\right\vert . \label{ip10}%
\end{equation}
For every $x,y\in\mathbb{R}^{d}$ with $\left\vert x-y\right\vert \leq1$\
\begin{align}
\left\vert f(y)-f(x)\right\vert  &  \leq\frac{1}{\Phi(x)}(1+\frac
{\widehat{\nabla}\Phi(x)}{\Phi(x)})\times\left\Vert \Phi f\right\Vert
_{1,p}\times\left\vert x-y\right\vert ^{\chi}\label{ip11}\\
&  \leq\frac{1}{\Phi(x)}(1+\frac{\widehat{\nabla}\Phi(x)}{\Phi(x)}%
)\times\Gamma_{\Phi,\theta}(1,h,p)\times\left\vert x-y\right\vert ^{\chi}.
\label{ip11a}%
\end{align}

\end{theorem}

\begin{remark}
The constants in the previous theorem depend on $\Phi$ by means of the
constants $c_{1}$ and $c_{2}$ which appear in the property (\ref{ip4}).
\end{remark}

\begin{remark}
The estimates in the point \textbf{C} in the above theorem are quite precise
and this is important in order to prove Theorem \ref{Bound}. But, roughly
speaking, (\ref{ip9}) reads
\[
\left\vert f(x)\right\vert \leq\frac{C}{\Phi(x)}\times(E(\Phi^{\theta
}(F))+\sup_{\varepsilon>0}E(\Phi^{\theta}(F_{\varepsilon}))).
\]
The constant $C$ depends on $h_{\ast},\widehat{H}_{2h_{\ast}+1+d,p_{\ast}%
}^{1/2h},\beta,C_{\ast},d$ and $p.$ This version is less precise but focus on
the following basic fact: if one is able to estimate the moments
$E(\Phi^{\theta}(F))$ and $E(\Phi^{\theta}(F_{\varepsilon}))$ then one obtains
the upper bound of $f$ by $\Phi^{-1}(x).$ This means that one is able to
translate moment estimates in terms of upper bounds for the density function.
\end{remark}

Let us try to give the heuristic which is behind the above criterion. Suppose
for simplicity that we are in dimension $d=1$ and that $F$ itself satisfies
the integration by parts formula (\ref{ip1}) (with $G=1).$ We formally write
\[
f_{F}(x)=E(\delta_{0}(F-x))=E(1_{[0,\infty)}^{\prime}(F-x))=E(1_{[0,\infty
)}(F-x)H_{1}).
\]
Using regularization by convolution the above reasoning may be done rigorously
and it proves that $P(F\in dx)=f_{F}(x)dx.$ Let us now compute the upper
bounds. Take for example $\Phi(x)=e^{\left\vert x\right\vert }.$ Then, using
Schwarz's inequality first and Chebishev's inequality then%
\[
f_{F}(x)\leq P^{1/2}(F\geq x)\left\Vert H_{1}\right\Vert _{2}=P^{1/2}%
(e^{\left\vert F\right\vert }\geq e^{\left\vert x\right\vert })\left\Vert
H_{1}\right\Vert _{2}\leq e^{-\frac{1}{2}\left\vert x\right\vert
}E(e^{\left\vert F\right\vert })^{1/2}\left\Vert H_{1}\right\Vert _{2}%
\]
so we obtain an estimate of type (\ref{ip9}).

The classical probabilistic way to obtain integration by parts formulas of
type (\ref{ip1}) is to use Malliavin calculus - our approach is strongly
inspired from this methodology, but however, at a certain point, takes a
completely different direction. Malliavin calculus is an infinite differential
calculus settled in the following way. One considers a class of "simple
functionals" which are "finite dimensional" objects. For them one defines a
derivative operator $D$ and a divergence operator $L$ using the classical
differential operators in finite dimension. Then one defines the extension of
these operators in infinite dimension: a general functional $F$ is in the
domain of $D$ (respectively of $L)$ if one may find a sequence of simple
functionals $F_{\varepsilon}$ such that $F_{\varepsilon}\rightarrow F$ in
$L^{2}$ and $DF_{\varepsilon}\rightarrow U$ in $L^{2}$ (respectively
$LF_{\varepsilon}\rightarrow V).$ Then one defines $DF=\lim_{\varepsilon
\rightarrow0}DF_{\varepsilon}$ respectively $LF=\lim_{\varepsilon\rightarrow
0}LF_{\varepsilon}.$ These operators are used in order to built the weight
$H_{1}$ in the integration by parts formula $E(\varphi^{\prime}(F))=E(\varphi
(F)H_{1}).$ In our approach we also settle a finite dimensional calculus as
above so we define $DF_{\varepsilon}$ and $LF_{\varepsilon}$ for a finite
dimensional $F_{\varepsilon}$ and we use them in order to obtain
$E(\varphi^{\prime}(F_{\varepsilon}))=E(\varphi(F_{\varepsilon}%
)H_{1,\varepsilon})$ (see Section 5). But in our framework $LF_{\varepsilon
}\uparrow\infty$ (see Remark \ref{D} for more details) so $F$ is no more in
the domain of $L.$ So the second step in Malliavin's methodology brakes down:
there is no infinite dimensional calculus available here. And we have
$H_{1,\varepsilon}\uparrow\infty.$ But, if we are able to obtain the estimates
given in (\ref{ip2}) and in (\ref{ip3}) and if the equilibrium condition
(\ref{ip6}) is verified, then we still obtain the regularity of the law of $F$
and the upper bound for its density. This is the object of Theorem
\ref{Criterion} . The ideas of this criterion origin in \cite{BCpreprint} and
the proof is given in Section 4.

In the years $80^{\prime}th$ starting with the papers \cite{[Bi]}, \cite{[L]}
and \cite{[BGJ]} a version of Malliavin calculus for Poisson point measures
has been developed and successfully used in order to study the regularity of
the solutions of $SDE^{\prime}s$ with jumps (see also \cite{[I]} and
\cite{[BD]} and the references there for recent developments in this area). In
the above papers the extension form the finite dimensional calculus to the
infinite dimensional one is successfully done and the limit $\lim
_{\varepsilon\rightarrow0}LF_{\varepsilon}=LF$ exists. Although the finite
dimensional calculus developed in our paper is similar, in our framework
$LF_{\varepsilon}$ blows up as $\varepsilon\rightarrow0.$ This is because the
law of the jumps in Boltzmann equation (and more generally in the framework of
Piecewise Deterministic Markov Processes) depends on the position of the
particle before the jump while for usual $SDE^{\prime}s,$ the law of the jumps
is independent of the position of the particle (see Section 5 for details).

The proof of Theorem \ref{L} is based on the criterion given in Theorem
\ref{Criterion}. In order to do it, following Tanaka \cite{T},\ we introduce a
stochastic equation which represents the probabilistic representation of the
Boltzmann equation and we construct some regularized version of this equation.
The solutions of these equations play the role of $F$ and of $F_{\varepsilon}$
in Theorem \ref{Criterion}. Then we recall two results from \cite{[BF]}: the
first one permits to estimate the error in (\ref{ip3}) and the second one
gives the integration by parts formula (\ref{ip1}) and the estimates in
(\ref{ip2}).

The paper is organized as follows: in Section 2 we recall the results from
\cite{[BF]} and in Section 3 we prove Theorem \ref{L}, Corollary \ref{R}
and\ Theorem \ref{Bound}\ (starting from the general criterion given in
Theorem \ref{Criterion}). In Section 4 we prove Theorem \ref{Criterion}. In
Section 5 we give an overview of the results from \cite{[BF]} and we precise
the changes which are necessary in order to obtain an explicit expression for
the dependence with respect to $t$ of the constants in the main estimates.

\bigskip

\section{Preliminary results}

In this section we present some results from \cite{[BF]}. Throughout this
section we fix $\nu\in\lbrack0,\frac{1}{2}),\gamma\in\lbrack0,1]$ and
$\lambda\in(\gamma,2)$ and the corresponding solution $f_{t}(dv)$ of
(\ref{bo1}) (which exists and is unique). In \cite{[BF]} (following the ideas
from \cite{T}) one gives the probabilistic interpretation of the equation
(\ref{bo1}). We recall this now. Let $E=[-\frac{\pi}{2},\frac{\pi}{2}%
]\times\mathbb{R}^{2}$ and let $N(dt,d\theta,dv,du)$ be a Poisson point
measure on $E\times\mathbb{R}_{+}$ with intensity measure $dt\times
b(\theta)d\theta\times f_{t}(dv)\times du.$ Consider also the matrix
\[
A(\theta)=\frac{1}{2}\left(
\begin{tabular}
[c]{ll}%
$\cos\theta-1$ & $-\sin\theta$\\
$\sin\theta$ & $\cos\theta-1$%
\end{tabular}
\ \ \right)  =\frac{1}{2}(R_{\theta}-I).
\]
We are interested in the equation
\begin{equation}
V_{t}=V_{0}+\int_{0}^{t}\int_{E\times\mathbb{R}_{+}}A(\theta)(V_{s-}%
-v)1_{\{u\leq\left\vert V_{s-}-v\right\vert ^{\gamma}\}}N(ds,d\theta,dv,du)
\label{bo3}%
\end{equation}
with $P(V_{0}\in dv)=f_{0}(dv).$ Proposition 2.1 in \cite{[BF]} asserts that
the equation (\ref{bo3}) has a unique c\`{a}dl\`{a}g solution $(V_{t}%
)_{t\geq0}$ and $P(V_{t}\in dv)=f_{t}(dv)$ (in this sense $V_{t}$ represents
the probabilistic representation for $f_{t}).$

In order to handle the equation (\ref{bo3})\ we face several difficulties: the
derivatives of the function $w\rightarrow\left\vert w-v\right\vert ^{\gamma}%
$\ blow up in the neighborhood of $v$ - so we have to use a regularization
procedure. Moreover, this function is unbounded and so we use a truncation
argument. Finally, the measure $\theta^{-(1+\nu)}d\theta$ has infinite mass,
and it is convenient to use a truncation argument also. We follow here the
ideas and results from \cite{[BF]}. We fix%
\begin{equation}
\eta_{0}\in(\frac{1}{\lambda},\frac{1}{\gamma\vee\nu})\quad and\quad
\Gamma_{\varepsilon}=(\ln\frac{1}{\varepsilon})^{\eta_{0}}. \label{bo3c}%
\end{equation}
Since $\gamma\eta_{0}>1$ we have, for every $C\geq1$ and $a>0$
\begin{equation}
\overline{\lim_{\varepsilon\rightarrow0}}\varepsilon^{a}e^{C\Gamma
_{\varepsilon}^{\gamma}}=0. \label{bo3a}%
\end{equation}
So $e^{C\Gamma_{\varepsilon}^{\gamma}}\leq\varepsilon^{-a}$ for sufficiently
small $\varepsilon.$ Moreover, if $\kappa>0$ is such that $\kappa\eta_{0}>1,$
then for every $A\geq1$%
\begin{equation}
\overline{\lim_{\varepsilon\rightarrow0}}\varepsilon^{-A}e^{-\Gamma
_{\varepsilon}^{\kappa}}=0. \label{bo3b}%
\end{equation}
So $e^{-\Gamma_{\varepsilon}^{\kappa}}\leq\varepsilon^{A}$ for sufficiently
small $\varepsilon.$

We construct the following approximation. We consider a $C^{\infty}$ even non
negative function $\chi$ supported by $[-1,1]$ and such that $\int_{R}%
\chi(x)dx=1$\ \ and we define%
\begin{equation}
\varphi_{\varepsilon}(x)=\int_{R}((y\vee2\varepsilon)\wedge\Gamma
_{\varepsilon})\frac{\chi((x-y)/\varepsilon)}{\varepsilon}dy. \label{bo5}%
\end{equation}
Observe that we have $2\varepsilon\leq\varphi_{\varepsilon}(x)\leq
\Gamma_{\varepsilon}$ for every $x\in\mathbb{R},\varphi_{\varepsilon}(x)=x$
for $x\in(3\varepsilon,\Gamma_{\varepsilon}-1),\varphi_{\varepsilon
}(x)=2\varepsilon$ for $x\in(0,\varepsilon)$ and $\varphi_{\varepsilon
}(x)=\Gamma_{\varepsilon}$ for $x\in(\Gamma_{\varepsilon},\infty).$ To the cut
off function $\varphi_{\varepsilon}$ one associates the equation
\begin{equation}
V_{t}^{\varepsilon}=V_{0}+\int_{0}^{t}\int_{E\times\mathbb{R}_{+}}%
A(\theta)(V_{s-}^{\varepsilon}-v)1_{\{u\leq\varphi_{\varepsilon}^{\gamma
}(\left\vert V_{s-}^{\varepsilon}-v\right\vert )\}}N(ds,d\theta,dv,du).
\label{bo7}%
\end{equation}
We construct a second approximation: for $\zeta>0$ we consider a smooth
cut-off function $I_{\zeta}$\ which is a smooth version of $1_{\{\left\vert
\theta\right\vert >\zeta\}}$ (the precise definition is given in (\ref{a2}))
and we associate the equation%
\begin{equation}
V_{t}^{\varepsilon,\zeta}=V_{0}+\int_{0}^{t}\int_{E\times\mathbb{R}_{+}%
}A(\theta)(V_{s-}^{\varepsilon,\zeta}-v)1_{\{u\leq\varphi_{\varepsilon
}^{\gamma}(\left\vert V_{s-}^{\varepsilon,\zeta}-v\right\vert )\}}I_{\zeta
}(\theta)N(ds,d\theta,dv,du). \label{bo7'}%
\end{equation}
We state now a property which will be used in the following: given $\alpha
\in\lbrack0,2]$ and $\kappa\geq0$ \ there exists $K\geq1$ such that for every
$w\in\mathbb{R}^{2},t_{0}>0$ and every $0<\varepsilon<1$%
\begin{equation}
(A_{\alpha,\kappa})\quad\sup_{t_{0}\leq t\leq T}f_{t}(Ball(w,\varepsilon
))\leq\frac{K}{t_{0}^{\kappa}}\varepsilon^{\alpha}. \label{bo8}%
\end{equation}
Since $f_{t}(dv)$ is a probability measure, this property is always verified
with $K=1,\alpha=0$ and $\kappa=0.$

In Proposition 2.1 from \cite{[BF]} one proves that the equations (\ref{bo7})
and (\ref{bo7'}) have a unique solution and
\begin{equation}
E\left\vert V_{t}^{\varepsilon,\zeta}-V_{t}^{\varepsilon}\right\vert \leq
C_{T}e^{C\Gamma_{\varepsilon}^{\gamma}}\times\zeta^{1-\nu}\times t\quad\forall
t\leq T. \label{bo8a}%
\end{equation}
Moreover, if $(A_{\alpha,\kappa})$ holds, then
\begin{equation}
E\left\vert V_{t}-V_{t}^{\varepsilon}\right\vert \leq C_{T}e^{C\Gamma
_{\varepsilon}^{\gamma}}\times\varepsilon^{1+\gamma+\alpha}\times t^{1-\kappa
}\quad\forall t\leq T. \label{bo8b}%
\end{equation}
We stress that in \cite{[BF]} the explicit dependence on the time $t$ does not
appear in the right hand side of the above estimates - but a quick glance to
the proof shows that we have the dependence on $t$ as in (\ref{bo8a}) and in
(\ref{bo8b}) ( this is important if we look to short time behavior). Moreover,
in the same proposition one proves that for every $0<\lambda^{\prime}<\lambda$
there exists some $\varepsilon_{0}>0$ such that%
\begin{equation}
\sup_{\varepsilon\leq\varepsilon_{0}}\sup_{\zeta\leq1}E(\sup_{t\leq
T}(e^{\left\vert V_{t}\right\vert ^{\lambda^{\prime}}}+e^{\left\vert
V_{t}^{\varepsilon}\right\vert ^{\lambda^{\prime}}}+e^{\left\vert V_{t}%
^{\zeta,\varepsilon}\right\vert ^{\lambda^{\prime}}})=:C(\lambda^{\prime
})<\infty. \label{bo4}%
\end{equation}

Finally in Theorem 4.1 in \cite{[BF]} one proves an integration by parts
formula that we present now. One defines (see (4.1) and (4.2) in \cite{[BF]})
a random process $G_{t}^{\varepsilon,\zeta}$ which verifies
\begin{equation}
1_{\{\sup_{s\leq t}\left\vert V_{s}^{\zeta,\varepsilon}\right\vert \leq
\Gamma_{\varepsilon}-1\}}\leq G_{t}^{\varepsilon,\zeta}\leq1_{\{\sup_{s\leq
t}\left\vert V_{s}^{\zeta,\varepsilon}\right\vert \leq\Gamma_{\varepsilon}\}}.
\label{bo8c}%
\end{equation}
The precise form of $G_{t}^{\varepsilon,\zeta}$ is not important here - the
only property which we need is (\ref{bo8c}). Moreover, since the law of
$V_{t}^{\varepsilon,\zeta}$ is not absolutely continuous we use the following
regularization procedure. One considers a two dimensional standard normal
random variable $Z$ and denotes
\begin{equation}
F_{t}^{\varepsilon,\zeta}=\sqrt{u_{\zeta}(t)}Z+V_{t}^{\varepsilon,\zeta
}\ \quad with\quad u_{\zeta}(t)=t\zeta^{4+\nu}. \label{bo8d}%
\end{equation}
Then one proves (see (4.3) and (4.4) in \cite{[BF]}) that for every
multi-index $\beta\in\{1,2\}^{q}$ there exists a random variable $K_{\beta
}(F_{t}^{\varepsilon,\zeta},G_{t}^{\varepsilon,\zeta})$ such that for every
function $\psi\in C^{q}(\mathbb{R}^{2})$%
\begin{equation}
E(\partial^{\beta}\psi(F_{t}^{\varepsilon,\zeta})G_{t}^{\varepsilon,\zeta
})=E(\psi(F_{t}^{\varepsilon,\zeta})K_{\beta}(F_{t}^{\varepsilon,\zeta}%
,G_{t}^{\varepsilon,\zeta})). \label{bo9a}%
\end{equation}
One also proves that for every $q\in N$ and every $\kappa\in(\frac{1}{\eta
_{0}},\lambda)$ one may find a constant $C$ (depending on $q$ and $\kappa$
only) such that for every $p\geq1$
\begin{equation}
\left\Vert K_{\beta}(F_{t}^{\varepsilon,\zeta},G_{t}^{\varepsilon,\zeta
})\right\Vert _{p}\leq\frac{C}{t^{\frac{2+\nu}{\nu}(12q-4)}}e^{C\Gamma
_{\varepsilon}^{\gamma}}(\varepsilon^{-q}\zeta^{-\nu q}+e^{-\Gamma
_{\varepsilon}^{\kappa}}\zeta^{-2\nu q}). \label{bo9b}%
\end{equation}
In particular this gives for every function $\psi\in C^{q}(\mathbb{R}^{2})$
and every multi-index $\beta\in\{1,2\}^{q}$%
\begin{equation}
\left\vert E(\partial^{\beta}\psi(F_{t}^{\varepsilon,\zeta})G_{t}%
^{\varepsilon,\zeta})\right\vert \leq\frac{C}{t^{\frac{2+\nu}{\nu}(12q-4)}%
}e^{C\Gamma_{\varepsilon}^{\gamma}}(\varepsilon^{-q}\zeta^{-\nu q}%
+e^{-\Gamma_{\varepsilon}^{\kappa}}\zeta^{-2\nu q})\times\left\Vert
\psi\right\Vert _{\infty}. \label{bo9}%
\end{equation}

The proof of (\ref{bo9a}) and (\ref{bo9b}) is based on a Malliavin type
calculus for jump processes and is quit technical. In Section 5 we give an
overview of the objects which come on in this proof and on the main estimates
which are needed. In particular we mention that in Theorem 4.1 in \cite{[BF]}
the dependence with respect to $t$ in the right hand side of (\ref{bo9b}) is
not explicit - at the end of Section 5 we precise this dependence (see
(\ref{a10})).

\bigskip

\section{Proof of Theorem \ref{L}}

In the following we adapt the results presented in the previous section to our
specific goals. Recall that the parameters $\nu\in\lbrack0,\frac{1}{2}%
),\gamma\in\lbrack0,1]$ and $\lambda\in(\gamma,2)$ are given and characterizes
the unique solution $f_{t}(dv)$ of (\ref{bo1}). Suppose that $(A_{\alpha
,\kappa})$\ (see (\ref{bo8})) holds for some $\alpha\geq0$ and $\kappa\geq0$.
In order to equilibrate the errors in (\ref{bo8a}) and (\ref{bo8b}) we take%
\[
\zeta=\zeta_{\alpha}(\varepsilon)=\varepsilon^{(1+\gamma+\alpha)/(1-\nu)}.
\]
With this choice, for every $c>0,$ we may find $C\geq1$ (depending on $c,$ see
((\ref{bo3a})) such that%
\begin{equation}
E(\left\vert V_{t}-V_{t}^{\varepsilon,\zeta_{\alpha}(\varepsilon)}\right\vert
)\leq\frac{C}{t^{\kappa-1}}\times e^{C\Gamma_{\varepsilon}^{\gamma}}%
\times\varepsilon^{1+\gamma+\alpha}\leq\frac{C}{t^{\kappa-1}}\varepsilon
^{1+\gamma+\alpha-c}. \label{cc7a}%
\end{equation}

Recall that $\eta_{0}$ is given in (\ref{bo3c}) and$\ \lambda\eta_{0}>1.$ So
we may choose (and fix) some $\lambda^{\prime}\in(\frac{1}{\eta_{0}}%
,\lambda).$ We work with the function $\Phi_{\lambda^{\prime}}$ given in
(\ref{Lp2'}) and we define
\[
g_{t}(dv)=\Phi_{\lambda^{\prime}}(v)f_{t}(dv).
\]
Moreover, for $\varepsilon>0,$ we recall that $F_{t}^{\varepsilon,\zeta}$ and
$G_{t}^{\varepsilon,\zeta}$ are given in (\ref{bo8d}) and (\ref{bo8c}), and we
define $f_{t}^{\varepsilon,\alpha}(dv)$ and $g_{t}^{\varepsilon,\alpha}(dv)$
by%
\[
\int\psi(v)f_{t}^{\varepsilon,\alpha}(dv)=E(\psi(F_{t}^{\varepsilon
,\zeta_{\alpha}(\varepsilon)})G_{t}^{\varepsilon,\zeta_{\alpha}(\varepsilon
)}),\quad g_{t}^{\varepsilon,\alpha}(dv)=\Phi_{\lambda^{\prime}}%
(v)f_{t}^{\varepsilon,\alpha}(dv).
\]
In (\ref{Lp1}) we introduced the function $\varphi.$ Notice that $\varphi$
solves the equation%
\begin{equation}
1+\gamma+\alpha-(1+\varphi(\alpha))\frac{1+\nu(\gamma+\alpha)}{1-\nu}=0.
\label{cc8a}%
\end{equation}
We construct the sequences
\begin{equation}
\alpha_{k+1}=\varphi(\alpha_{k}),\quad\kappa_{k+1}=\kappa_{k}-1+\frac
{13(2+\nu)}{\nu}(1+\alpha_{k+1}) \label{Lp2}%
\end{equation}
with $\alpha_{0}=0$ and $\kappa_{0}=0.$ Direct computation shows that
$\varphi^{\prime}(\alpha)>0$ for every $\alpha,$ so $\varphi$ is strictly
increasing. We will assume in the following that $\alpha_{1}=\varphi
(0)>0=\alpha_{0}$ and this implies $\alpha_{k+1}>\alpha_{k}$ for every $k.$ It
follows that $\alpha_{k}\uparrow\alpha_{\ast}$ solution of $\varphi
(\alpha_{\ast})=\alpha_{\ast}$ (see (\ref{alpha}) for the explicit value of
$\alpha_{\ast}).$ Notice also that $\alpha_{1}=\varphi(0)>0$ is equivalent
with $\nu<\frac{\gamma}{2\gamma+1}$ and $\varphi(2)>2$ is equivalent with
$\nu<\frac{\gamma}{4\gamma+9}.$

We know that $(A_{0,0})$ holds. Our aim now is to employ Theorem
\ref{Criterion} in order to obtain $(A_{\alpha,\kappa})$ for $\alpha$ as large
as possible.

\begin{lemma}
\label{A}\textbf{A}. Let $q\in\mathbb{N},\alpha\in\lbrack0,2]$ and $\kappa
\geq0$ be given. Suppose that $(A_{\alpha,\kappa})$ holds with $\varphi
(\alpha)>q,$ and take $p>1$ such that
\begin{equation}
q+\frac{2}{p_{\ast}}<\varphi(\alpha). \label{cc5}%
\end{equation}
Then $f_{t}(dv)=f_{t}(v)dv$ with $f_{t}\in W^{q,p}.$ Moreover, for every
$\lambda^{\prime}<\lambda$ there exists $C\geq1$ such that
\begin{equation}
\left\Vert \Phi_{\lambda^{\prime}}f_{t}\right\Vert _{q,p}\leq\frac
{C}{t^{\kappa-1+\frac{13(2+\nu)}{\nu}(1+\varphi(\alpha))}}. \label{cc6}%
\end{equation}
\textbf{B} Suppose that $(A_{\alpha,\kappa})$ holds and $\varphi(\alpha)>0.$
Then $(A_{\alpha^{\prime},\kappa^{\prime}})$ holds for every $\alpha^{\prime
}<\varphi(\alpha)\wedge2$ with
\begin{equation}
\kappa^{\prime}=\kappa-1+\frac{13(2+\nu)}{\nu}(1+\varphi(\alpha)). \label{cc7}%
\end{equation}
\textbf{C}. Let $\alpha_{k},\kappa_{k},k\in\mathbb{N}$ be the sequences
defined in (\ref{Lp2}). Suppose that $\varphi(0)>0.$ Then, for each
$k\in\mathbb{N}_{\ast}$ the property $(A_{\alpha,\kappa_{k}})$ holds for every
$\alpha<\alpha_{k}\wedge2.$

\textbf{D}. Suppose that $\varphi(0)>0.$ Let $k,q\in\mathbb{N}$ and $p>1$ be
such that
\begin{equation}
q+\frac{2}{p_{\ast}}<\varphi(\alpha_{k}\wedge2)=\alpha_{k+1}\wedge\varphi(2).
\label{cc7b}%
\end{equation}
Then
\begin{equation}
\left\Vert \Phi_{\lambda^{\prime}}f_{t}\right\Vert _{q,p}\leq\frac
{C}{t^{\kappa_{k+1}}}. \label{cc7c}%
\end{equation}

\end{lemma}

\textbf{Proof of A.} We will use Theorem \ref{Criterion} with $d=2,$ and
$F_{\varepsilon}=F_{t}^{\varepsilon,\zeta_{\alpha}(\varepsilon)}%
,G_{\varepsilon}=G_{t}^{\varepsilon,\zeta_{\alpha}(\varepsilon)}.$ So we
verify the hypothesis there.

\textbf{Step 1}. By (\ref{bo9a}) we know that the integration by parts formula
(\ref{ip1}) holds with $H_{\beta,\varepsilon}=K_{\beta}(F_{t}^{\varepsilon
,\zeta_{\alpha}(\varepsilon)},G_{t}^{\varepsilon,\zeta_{\alpha}(\varepsilon
)}).$ By (\ref{bo9b}) we obtain for every $\kappa\in(\frac{1}{\eta_{0}%
},\lambda)$ (with $\zeta=\zeta_{\alpha}(\varepsilon))$
\[
\sup_{\left\vert \beta\right\vert \leq q}\left\Vert H_{\beta,\varepsilon
}\right\Vert _{p}\leq\frac{C}{t^{\frac{2+\nu}{\nu}(12q-4)}}e^{C\Gamma
_{\varepsilon}^{\gamma}}(\varepsilon^{-q}\zeta^{-\nu q}+e^{-\Gamma
_{\varepsilon}^{\kappa}}\zeta^{-2\nu q}).
\]
We use (\ref{bo3a}) and (\ref{bo3b}) in order to obtain%
\[
e^{C\Gamma_{\varepsilon}^{\gamma}}(\varepsilon^{-q}\zeta^{-\nu q}%
+e^{-\Gamma_{\varepsilon}^{\kappa}}\zeta^{-2\nu q})\leq C\varepsilon
^{-c}((\varepsilon\zeta^{\nu})^{-q}+\varepsilon^{A}\varepsilon^{-2\nu
q(1+\gamma+\alpha)/(1-\nu)})
\]
for every $c>0$ and $A\geq1.$ Notice that $\varepsilon\zeta_{\alpha}^{\nu
}(\varepsilon)=\varepsilon^{(1+\nu(q+\alpha))/(1-\nu)}.$ Taking $A\geq2\nu
q(1+\gamma+\alpha)/(1-\nu)$ we obtain%

\[
\sup_{\left\vert \beta\right\vert \leq q}\left\Vert H_{\beta,\varepsilon
}\right\Vert _{p}\leq\frac{C}{t^{\frac{2+\nu}{\nu}(12q-4)}}\times
\varepsilon^{-q\times\frac{1+\nu(\gamma+\alpha)}{1-\nu}-c}%
\]
and this means that (\ref{ip2}) is verified with
\[
\widehat{H}_{q,p}=\frac{C}{t^{\frac{2+\nu}{\nu}(12q-4)}},\quad b=\frac
{1+\nu(\gamma+\alpha)}{1-\nu},\quad a=\frac{c}{b}.
\]
Let $\delta>0.$ Taking $h$ sufficiently large we have $\frac{1}{2h}%
(12(2h+q+2)-4)\leq12(1+\delta)$ so that%
\begin{equation}
\widehat{H}_{2h+q+2,p}^{1/2h}\leq\frac{C}{t^{\frac{2+\nu}{\nu}\times
12(1+\delta)}}. \label{c7a}%
\end{equation}

\textbf{Step 2}. Let us verify (\ref{ip3}). Using (\ref{bo8c}) and
(\ref{bo4})
\begin{align*}
\left\Vert 1-G_{t}^{\varepsilon,\zeta_{\alpha}(\varepsilon)}\right\Vert _{2}
&  \leq P^{1/2}(\sup_{s\leq t}\left\vert V_{s}^{\varepsilon,\zeta_{\alpha
}(\varepsilon)}\right\vert \geq\Gamma_{\varepsilon})\\
&  \leq Ce^{-\frac{1}{2}\Gamma_{\varepsilon}^{\lambda^{\prime}}}(E(\sup_{s\leq
t}e^{\left\vert V_{s}^{\varepsilon,\zeta_{\alpha}(\varepsilon)}\right\vert
^{\lambda^{\prime}}}))^{1/2}\leq Ce^{-\frac{1}{2}\Gamma_{\varepsilon}%
^{\lambda^{\prime}}}\leq C\varepsilon^{A}.
\end{align*}
The last inequality is true for any $A\geq1.$ It is a consequence of
$\lambda^{\prime}\eta_{0}>1$ and of (\ref{bo3b}).

Recall that $F_{t}^{\varepsilon,\zeta_{\alpha}(\varepsilon)}$ is defined in
(\ref{bo8d}). We have
\[
E(\left\vert V_{t}^{\varepsilon,\zeta_{\alpha}(\varepsilon)}-F_{t}%
^{\varepsilon,\zeta_{\alpha}(\varepsilon)}\right\vert )\leq C\zeta_{\alpha
}^{\frac{4+\nu}{2}}(\varepsilon)=C\varepsilon^{(1+\gamma+\alpha)\frac{4+\nu
}{2(1-\nu)}}\leq C\varepsilon^{1+\gamma+\alpha}.
\]
Then, as a consequence of (\ref{cc7a}), for every $c>0$ we obtain%
\[
E(\left\vert V_{t}-F_{t}^{\varepsilon,\zeta_{\alpha}(\varepsilon)}\right\vert
)\leq\frac{C}{t^{\kappa-1}}\varepsilon^{1+\gamma+\alpha-c}.
\]
We conclude that (\ref{ip3}) holds with%
\begin{equation}
C_{\ast}=\frac{C}{t^{\kappa-1}},\quad\beta=1+\gamma+\alpha-c. \label{c7b}%
\end{equation}

\textbf{Step 3}. Now (\ref{cc5}) ensures that, for sufficiently small $c>0,$
\begin{align*}
\beta-b(1+q+\frac{2}{p_{\ast}})  &  =1+\gamma+\alpha-c-(1+q+\frac{2}{p_{\ast}%
})\frac{1+\nu(\gamma+\alpha)}{1-\nu}\\
&  >1+\gamma+\alpha-(1+\varphi(\alpha))\frac{1+\nu(\gamma+\alpha)}{1-\nu}=0
\end{align*}
the last equality being a consequence of (\ref{cc8a}) (this is the motivation
of choosing $\varphi$ to be the solution of this equation).

So (\ref{ip6})\ holds (with $d=2)$ and we are able to use Theorem
\ref{Criterion}. \ Notice that for every $\theta\geq1$ and every
$\lambda^{\prime}<\lambda^{\prime\prime}$ one may find $C$ such that
$\Phi_{\lambda^{\prime}}^{\theta}\leq C\Phi_{\lambda^{\prime\prime}}.$ So
(\ref{bo4}) gives $C_{\theta}(\Phi)<\infty$ (see (\ref{ip5}) for the
definition of $C_{\theta}(\Phi)).$ By (\ref{ip7})%
\[
\left\Vert \Phi_{\lambda^{\prime}}f_{t}\right\Vert _{q,p}\leq C\times(C_{\ast
}+C_{\theta}(\Phi))\times\left(  h^{2b}\widehat{H}_{2h+q+d,p_{\ast}}%
^{1/2h}\right)  ^{(1+\delta)(1+q+2/p_{\ast})}.
\]
We denote $C^{\prime}=C_{\theta}(\Phi)\times h^{2b(1+\delta)(1+q+2/p_{\ast}%
)}.$ Then, using (\ref{c7b}) and (\ref{c7a})%
\begin{align*}
\left\Vert \Phi_{\lambda^{\prime}}f_{t}\right\Vert _{q,p}  &  \leq\frac
{C}{t^{\kappa-1+\frac{2+\nu}{\nu}\times12(1+q+2/p_{\ast})(1+\delta)^{2}}}\\
&  \leq\frac{C}{t^{\kappa-1+\frac{2+\nu}{\nu}\times12(1+\varphi(\alpha
))(1+\delta)^{2}}}.
\end{align*}
We take $\delta>0$ sufficiently small so that $12(1+\delta)^{2}\leq13$ and
\textbf{A} is proved.

\textbf{Proof of B}.\ We use \textbf{A} with $q=0.$ Let $0<\alpha^{\prime
}<\varphi(\alpha)\wedge2.$ Since $\alpha^{\prime}<2$ we may find $p>1$ such
that $\frac{2}{p_{\ast}}=\alpha^{\prime}<\varphi(\alpha),$ so (\ref{cc5})
holds for this $p.$\ By(\ref{cc6}) $\left\Vert f_{t}\right\Vert _{p}%
\leq\left\Vert \Phi_{\lambda^{\prime}}f_{t}\right\Vert _{p}\leq Ct^{-\kappa
^{\prime}}$ with $\kappa^{\prime}$ given in (\ref{cc7}). Using H\"{o}lder's
inequality we get $(A_{\alpha^{\prime},\kappa^{\prime}}):$
\[
f_{t}(Ball(w,\varepsilon))\leq\left\Vert f_{t}\right\Vert _{p}\times
\varepsilon^{2/p_{\ast}}\leq Ct^{-\kappa^{\prime}}\times\varepsilon
^{\alpha^{\prime}}.
\]

\textbf{Proof of C. }Take first $k=1.$ We know that $(A_{0,0})$\ holds, and by
hypothesis, $\varphi(0)>0.$ Then, according to \textbf{B}, $(A_{\alpha
^{\prime},\kappa^{\prime}})$ holds for every $\alpha^{\prime}<\varphi
(0)\wedge2=\alpha_{1}\wedge2$ with $\kappa^{\prime}=0-1+\frac{13(2+\nu)}{\nu
}(1+\varphi(0))=\kappa_{1}.$ So our assertion holds for $k=1.$

Suppose now that the property is true for $k$ and let us check it for $k+1.$
Suppose first that $\alpha_{k}>2.$ Then $\alpha_{k+1}>\alpha_{k}>2.$ By the
recurrence hypothesis, for $\alpha<2=\alpha_{k}\wedge2=\alpha_{k+1}\wedge2$
the hypothesis $(A_{\alpha,\kappa_{k}})$ holds. Since $\kappa_{k+1}>\kappa
_{k},$the hypothesis $(A_{\alpha,\kappa_{k+1}})$ holds as well.

Suppose now that $\alpha_{k}\leq2$ and take $\alpha^{\prime}<\alpha
_{k+1}\wedge2=\varphi(\alpha_{k})\wedge2.$ Since $\varphi(\alpha
)\uparrow\varphi(\alpha_{k})$ as $\alpha\uparrow\alpha_{k},$ we may find
$\alpha<\alpha_{k}=\alpha_{k}\wedge2$ such that $\alpha^{\prime}%
<\varphi(\alpha)\wedge2.$ By the recurrence hypothesis we know that
$(A_{\alpha,\kappa_{k}})$ holds and then, using \textbf{B}, we obtain
$(A_{\alpha^{\prime},\kappa_{k+1}}).$

\textbf{D}. By (\ref{cc7b}) we may find $\alpha<\alpha_{k}\wedge2$ such that
$q+\frac{2}{p_{\ast}}<\varphi(\alpha).$ By \textbf{C} we know that
$(A_{\alpha,\kappa_{k}})$ holds. Then we may use \textbf{A} and (\ref{cc6})
gives (\ref{cc7c}). $\square$

\textbf{Proof of Theorem} \ref{L}. We will work with the sequences $\alpha
_{k}$ and $\kappa_{k}$ given in (\ref{Lp2}). Recall that $\alpha_{k}%
\uparrow\alpha_{\ast}$ with $\alpha_{\ast}=\varphi(\alpha_{\ast}).$ Direct
computation give
\begin{equation}
\alpha_{\ast}=\frac{-(\gamma+2)+\sqrt{(\gamma+2)^{2}+4(\frac{\gamma}{\nu
}-2\gamma-1)}}{2} \label{alpha}%
\end{equation}
and
\begin{align*}
\alpha_{\ast}  &  >0\quad\Leftrightarrow\quad\nu<\frac{\gamma}{2\gamma+1}\\
\alpha_{\ast}  &  >1\quad\Leftrightarrow\quad\nu<\frac{\gamma}{3\gamma+4}\\
\alpha_{\ast}  &  >2\quad\Leftrightarrow\quad\nu<\frac{\gamma}{4\gamma+9}.
\end{align*}
If $\nu<\frac{\gamma}{2\gamma+1}$ then $\varphi(0)>0$ so we may use the point
\textbf{A }in Lemma \ref{A} with $q=0.$ We obtain $f_{t}(dv)=f_{t}(v)dv$ so
the point \textbf{A} in Theorem \ref{L} is proved.

\textbf{Proof of} \textbf{B.b}. If $\frac{\gamma}{4\gamma+9}<\nu<\frac{\gamma
}{2\gamma+1}$ we have $\alpha_{\ast}\leq2$ so that $\alpha_{k}<2$ for every
$k\in N.$ If $p<\frac{2}{2-\alpha_{\ast}}$ then $\frac{2}{p_{\ast}}%
<\alpha_{\ast}$ so we may find $k$ such that $\frac{2}{p_{\ast}}<\alpha
_{k+1}<2.$ Using the point \textbf{D} in Lemma \ref{A} (with $q=0)$ we obtain
$\Phi_{\lambda^{\prime}}f_{t}\in L^{p}(R^{2}).$

\textbf{Proof of} \textbf{C.b}. If $\frac{\gamma}{4\gamma+9}<\nu<\frac{\gamma
}{3\gamma+4}$ we have $\alpha_{\ast}\in(1,2]$ so that $1<\frac{2}%
{3-\alpha_{\ast}}.$ We take $1<p<\frac{2}{3-\alpha_{\ast}}$ and then
$1+\frac{2}{p_{\ast}}<\alpha_{\ast}.$ We take $k$ sufficiently large in order
to have $1+\frac{2}{p_{\ast}}<\alpha_{k+1}$ and then, as above, by \textbf{D}
in Lemma \ref{A}, we obtain $\Phi_{\lambda^{\prime}}f_{t}\in W^{1,p}%
(\mathbb{R}^{2}).$

\textbf{Proof of} \textbf{B a }If $\nu<\frac{\gamma}{4\gamma+9}$ then
$\alpha_{\ast}>2.$ Recall that $\alpha_{k}\uparrow\alpha_{\ast}$ and define
$k_{\ast}=\min\{k:\alpha_{k}\geq2\}.$ By \textbf{C} in Lemma \ref{A}, for
every $\alpha<2$ the property $(A_{\alpha,\kappa_{k_{\ast}}})$ holds. We
denote $\psi(\alpha)=\varphi(\alpha)-\alpha.$ Notice that $\psi$ is decreasing
on $(0,2).$ For $k<k_{\ast}$ we have $\alpha_{k}<2=\alpha_{k_{\ast}}\wedge2$
so that $\psi(2)<\psi(\alpha_{k})=$ $\alpha_{k+1}-\alpha_{k}.$ It follows
that
\[
2>\alpha_{k_{\ast}-1}=\sum_{k=0}^{k_{\ast}-2}(\alpha_{k+1}-\alpha
_{k})>(k_{\ast}-1)\psi(2)
\]
which gives $k_{\ast}-1\leq2/\psi(2)$ and so $k_{\ast}+1\leq2(\varphi
(2)-1)/(\varphi(2)-2).$ Since $\alpha_{k}\leq\alpha_{\ast}=\varphi
(\alpha_{\ast})$ we have for every $k$ (see (\ref{Lp2}))
\[
\kappa_{k}\leq\kappa_{k-1}+\frac{13(2+\nu)}{\nu}(1+\varphi(\alpha_{\ast
}))-1\leq....\leq k(\frac{13(2+\nu)}{\nu}(1+\alpha_{\ast})-1).
\]
This yields
\[
\kappa_{k_{\ast}+1}\leq(k_{\ast}+1)\times(\frac{13(2+\nu)}{\nu}(1+\alpha
_{\ast})-1)\leq\frac{2(\varphi(2)-1)}{\varphi(2)-2}\left(  \frac
{13(1+\alpha_{\ast})(2+\nu)}{\nu}-1\right)  =\eta
\]
with $\eta$ given in (\ref{Lp3b}).

We use now the point \textbf{D} in Lemma \ref{A} with $q=0.$ Recall that
$\varphi(2)>2$ (see (\ref{Lp3d})) Since $\varphi(\alpha_{k_{\ast}}%
\wedge2)=\varphi(2)>2>\frac{2}{p_{\ast}}$ for every $p>1,$ we obtain $g_{t}%
\in\cap_{p>1}L^{p}(R^{2}).$ Moreover, taking $k=k_{\ast}$ in (\ref{cc7c}) we
get
\[
\left\Vert g_{t}\right\Vert _{p}\leq\frac{C}{t^{\kappa_{k_{\ast}+1}}}\leq
\frac{C}{t^{\eta}}.
\]
\textbf{Proof of C.a}. If $q\in\{1,2\}$, we need $\frac{2}{p_{\ast}}%
<\varphi(2)-q.$ This gives $p<2/(q+2-\varphi(2))=p_{q}$ with $p_{q},q=1,2$
given in (\ref{Lp4}). And using (\ref{cc7c}) we obtain%
\[
\left\Vert g_{t}\right\Vert _{q,p}\leq\frac{C}{t^{\kappa_{k_{\ast}+1}}}%
\leq\frac{C}{t^{\eta}}\quad for\quad p<p_{q}.
\]
$\square$

\textbf{Proof of Corollary} \ref{R} Recall that $\nu<\frac{\gamma}{2\gamma+1}$
is equivalent with $\varphi(0)>0.$ So we may find $p>1$ such that $\frac
{2}{p_{\ast}}<\varphi(0).$ Using \textbf{D} in Lemma \ref{A} with $q=0$ and
$k=0$ we obtain $\left\Vert \Phi_{\lambda^{\prime}}f_{t}\right\Vert _{p}%
\leq\frac{C}{t^{\kappa_{1}}}$ with $\kappa_{1}$ given in (\ref{Lp2})$\ $(which
coincides with $\kappa$ defined in (\ref{Lp6a})) Then%
\begin{align*}
f_{t}(B_{R}^{c}(0))  &  =\int1_{B_{R}^{c}(0)}(v)\Phi_{\lambda^{\prime}}%
^{-1}(v)\Phi_{\lambda^{\prime}}(v)f_{t}(v)dv\\
&  \leq(\int1_{B_{R}^{c}(0)}(v)e^{-p_{\ast}\left\vert v\right\vert
^{\lambda^{\prime}}}dv)^{1/p_{\ast}}\left\Vert \Phi_{\lambda^{\prime}}%
f_{t}\right\Vert _{p}\\
&  \leq e^{-\frac{1}{2}R^{\lambda^{\prime}}}(\int1_{B_{R}^{c}(0)}%
(v)e^{-\frac{p_{\ast}}{2}\left\vert v\right\vert ^{\lambda^{\prime}}%
}dv)^{1/p_{\ast}}\frac{C}{t^{\kappa}}\\
&  \leq\frac{C}{t^{\kappa}}e^{-\frac{1}{2}R^{\lambda^{\prime}}}.
\end{align*}
$\square$

\textbf{Proof of Theorem} \ref{Bound}. Since $\nu<\frac{\gamma}{4\gamma+9}$ we
may use \textbf{C.a} in Theorem \ref{L} with $q=1$ and we obtain $\left\Vert
\Phi_{\lambda^{\prime}}f\right\Vert _{1,p}\leq Ct^{-\eta}$ for every $p<p_{1}$
(with $p_{1}$ given in (\ref{Lp4}) and $\eta$ given in (\ref{Lp3b})). Notice
also that we have $p_{1}>2$ so we may use the point \textbf{C} in Theorem
\ref{Criterion}: (\ref{ip9}) gives (\ref{bo2a}). Notice also that if
$\lambda^{\prime}<\lambda^{\prime\prime}$ then $\widehat{\nabla}\Phi
_{\lambda^{\prime}}(x)\leq C\Phi_{\lambda^{\prime\prime}}(x).$ So (\ref{ip11})
gives (\ref{bo2bc}). $\square$

\section{Appendix: A regularity criterion based on interpolation}

\label{sect:3.1}

Let us first recall some results obtained in \cite{BCpreprint} concerning the
regularity of a measure $\mu$ on ${\mathbb{R}}^{d}$. Fix $k,q,h\in{\mathbb{N}%
}$, with $h\geq1$, and $p>1$ (we denote by $p_{\ast}$ the conjugate of $p)$.
For $f\in C^{\infty}(R^{d})$ we define
\begin{align}
\left\Vert f\right\Vert _{k,\infty}  &  =\sum_{0\leq\left\vert \alpha
\right\vert \leq k}\sup_{x\in R^{d}}\left\vert \partial^{\alpha}%
f(x)\right\vert ,\label{bo2bb}\\
\left\Vert f\right\Vert _{k,h,p}  &  =\sum_{0\leq\left\vert \alpha\right\vert
\leq k}(E(\int_{R^{d}}(1+\left\vert x\right\vert )^{h}\left\vert
\partial^{\alpha}f(x)\right\vert ^{p}dx))^{1/p}\label{bo2b}\\
\left\Vert f\right\Vert _{k,p}  &  =\left\Vert f\right\Vert _{k,0,p}%
=\sum_{0\leq\left\vert \alpha\right\vert \leq k}\left\Vert \partial^{\alpha
}f\right\Vert _{p}. \label{bo2bbb}%
\end{align}
Here $\alpha=(\alpha_{1},...,\alpha_{m})\in\{1,...,d\}^{m},$ is a multi-index
of length $\left\vert \alpha\right\vert =m$ and $\partial^{\alpha}$ is the
derivative associated to $\alpha.$ Moreover for two measures $\mu,\nu$ we
consider the distance%
\begin{equation}
d_{k}(\mu,\nu)=\sup\{\left\vert \int fd\mu-\int fd\nu\right\vert :\left\Vert
f\right\Vert _{k,\infty}\leq1\}. \label{bo2c}%
\end{equation}
For $k=0$ this is the total variation distance and for $k=1$ this is the
Fort\`{e}t Mourier distance.

For a finite measure $\mu$ and for a sequence of absolutely continuous finite
measures $\mu_{n}(dx)=f_{n}(x)dx$ with $f_{n}\in C^{2h+q}({\mathbb{R}}^{d}),$
we define%
\begin{equation}
\pi_{k,q,h,p}(\mu,(\mu_{n})_{n})=\sum_{n=0}^{\infty}2^{n(k+q+d/p_{\ast})}%
d_{k}(\mu,\mu_{n})+\sum_{n=0}^{\infty}\frac{1}{2^{2nh}}\left\Vert
f_{n}\right\Vert _{2h+q,2h,p} \label{reg3}%
\end{equation}
and%
\[
\overline{\pi}_{k,q,h,p}(\mu)=\inf\{\pi_{k,q,h,p}(\mu,(\mu_{n})_{n}):\mu
_{n}(dx)=f_{n}(x)dx,\quad f_{n}\in C^{2h+q}({\mathbb{R}}^{d})\}.
\]

\begin{remark}
Notice that $\pi_{k,q,h,p}$ is a particular case of $\pi_{k,q,h,\mathbf{e}}$
treated in \cite{BCpreprint}: just choose the Young function $\mathbf{e}%
(x)\equiv\mathbf{e}_{p}(x)=|x|^{p}$ (see Example 1 in \cite{BCpreprint}).
Moreover, $\pi_{k,q,h,p}$ is strongly related to interpolation spaces. More
precisely, $\overline{\pi}_{k,q,h,p}$ is equivalent with the interpolation
norm of order $\rho=\frac{k+q+d/p_{\ast}}{2h}$ between the spaces $W_{\ast
}^{k,\infty}$ (the dual of $W^{k,\infty})$ and $W^{2h+q,2h,p}=\{f:$
$\left\Vert f\right\Vert _{2h+q,2h,p}<\infty\}$ (see \cite{[BS]} for example).
This is proved in \cite{BCpreprint}, see Section 2.4 and Appendix B. So the
inequality (\ref{reg4}) below implies that the Sobolev space $W^{q,p}$ is
included in the above interpolation space. However we prefer to stick to an
elementary framework and to derive directly the consequences of (\ref{reg4}) -
see Lemma \ref{REG} and Lemma \ref{continuity} below.
\end{remark}

The following result is the key point in our approach:

\begin{lemma}
\label{continuity} Let $p>1,$ $k,q\in{\mathbb{N}}$ and $h\in\mathbb{N}_{\ast
}=\mathbb{N}\prime\{0\}$ be given. There exists a constant $C_{\ast}$
(depending on $k,q,h$ and $p$ only) such that the following holds. Let $\mu$
be a finite measure for which $\overline{\pi}_{k,q,h,p}(\mu)$ is finite. Then
$\mu(dx)=f(x)dx$ with $f\in W^{q,p}$ and
\begin{equation}
\left\Vert f\right\Vert _{q,p}\leq C_{\ast}\times\overline{\pi}_{k,q,h,p}%
(\mu). \label{reg4}%
\end{equation}

\end{lemma}

This is Proposition 2.5 in \cite{BCpreprint} in the particular case
$e(x)=e_{p}(x)=\left\vert x\right\vert ^{p}.$ See also Proposition 3.2.1 in
\cite{BCC} .\ So we will not give here the proof. We will use the following consequence:

\begin{lemma}
\label{REG} Let $p>1,$ $k,q\in{\mathbb{N}}$ and $h\in\mathbb{N}_{\ast}$ be
given and set
\begin{equation}
\rho_{h}(q):=\frac{k+q+d/p_{\ast}}{2h}. \label{reg5}%
\end{equation}
We consider an increasing sequence $\theta(n)\geq1,n\in{\mathbb{N}}$\ such
that $\lim_{n}\theta(n)=\infty$\ and $\theta(n+1)\leq\Theta\times\theta(n)$
for some constant $\Theta\geq1.$ Moreover, we consider a sequence of measures
$\mu_{n}(dx)=f_{n}(x)dx$ with $f_{n}\in C^{2h+q}({\mathbb{R}}^{d}%
),n\in{\mathbb{N}}$ such that
\begin{equation}
\left\Vert f_{n}\right\Vert _{2h+q,2h,p}\leq\theta(n). \label{reg9}%
\end{equation}
Let $\mu$ be a finite measure such that, for some $\eta>0,$%
\begin{equation}
\limsup_{n}d_{k}(\mu,\mu_{n})\times\theta^{\rho_{h}(q)+\eta}(n)<\infty.
\label{reg10}%
\end{equation}
Then $\mu(dx)=f(x)dx$ with $f\in W^{q,p}.$

Moreover, fix $n_{\ast}\in{\mathbb{N}}$, $\delta>0$ and $\eta>0$ such that
(\ref{reg10}) holds. We set
\begin{align}
A(\delta)  &  =\left\vert \mu\right\vert ({\mathbb{R}}^{d})\times
2^{l(\delta)(1+\delta)(q+k+d/p_{\ast})}\quad with\label{reg12'}\\
l(\delta)  &  =\min\{l\geq1:2^{l^{\prime}\times\frac{\delta}{1+\delta}}\geq
l^{\prime},\forall l^{\prime}\geq l\}\label{reg12'''}\\
B(\eta)  &  =\sum_{l=1}^{\infty}\frac{l^{2(q+k+d/p_{\ast}+\eta)}}{2^{2h\eta
l}},\label{reg12''}\\
C_{h,n_{\ast}}(\eta)  &  =\sup_{n\geq n_{\ast}}d_{k}(\mu,\mu_{n})\times
\theta^{\rho_{h}(q)+\eta}(n). \label{reg11}%
\end{align}
Then
\begin{equation}
\left\Vert f\right\Vert _{q,p}\leq C_{\ast}(\Theta+A(\delta)\theta(n_{\ast
})^{\rho_{h}(q)(1+\delta)}+B(\eta)C_{h,n_{\ast}}(\eta)) \label{reg12}%
\end{equation}
with $C_{\ast}$ the constant in (\ref{reg4}) and $\rho_{h}(q)$ given in
(\ref{reg5}).
\end{lemma}

\textbf{Proof of Lemma \ref{REG}}. We will produce a sequence of measures
$\nu_{l}(dx)=g_{l}(x)dx,l\in{\mathbb{N}}$ such that
\[
\pi_{k,q,h,p}(\mu,(\nu_{l})_{l})\leq\Theta+A(\delta)\theta(n_{\ast})^{\rho
_{h}(q)(1+\delta)}+B(\eta)C_{h,n_{\ast}}(\eta)<\infty.
\]
Then by Lemma \ref{continuity} one gets $\mu(dx)=f(x)dx$ with $f\in W^{q,p}, $
and (\ref{reg12}) follows from (\ref{reg4})). Let us stress that the $\nu_{l}%
$'s will be given by a suitable subsequence $\mu_{n(l)}$, $l\in{\mathbb{N}}$.

\textbf{Step 1}. We define
\[
n(l)=\min\{n:\theta(n)\geq\frac{2^{2hl}}{l^{2}}\}
\]
and we notice that
\begin{equation}
\frac{1}{\Theta}\theta(n(l))\leq\theta(n(l)-1)<\frac{2^{2hl}}{l^{2}}\leq
\theta(n(l)). \label{reg13}%
\end{equation}
Moreover we recall that $n_{\ast}$ is given and\ we define
\[
l_{\ast}=\min\{l:\frac{2^{2hl}}{l^{2}}\geq\theta(n_{\ast})\}.
\]
Since%
\[
\theta(n(l_{\ast}))\geq\frac{2^{2hl_{\ast}}}{l_{\ast}^{2}}\geq\theta(n_{\ast
})
\]
it follows that $n(l_{\ast})\geq n_{\ast}.$

We take now $\varepsilon(\delta)=\frac{h\delta}{1+\delta}$ which gives
$\frac{2h}{2(h-\varepsilon(\delta))}=1+\delta.$ And we take $l(\delta)\geq1$
such that $2^{l\delta/(1+\delta)}\geq l$ for $l\geq l(\delta)$ (see
(\ref{reg12'''}))$.$ Since $h\geq1$ it follows that $\varepsilon(\delta
)\geq\frac{\delta}{1+\delta}$ so that, for $l\geq l(\delta)$ we also have
$2^{l\varepsilon(\delta)}\geq l.$ Now we check that
\begin{equation}
2^{2(h-\varepsilon(\delta))l_{\ast}}\leq2^{2hl(\delta)}\theta(n_{\ast}).
\label{reg14}%
\end{equation}
If $l_{\ast}\leq l(\delta)$ then the inequality is evident (recall that
$\theta(n)\geq1$ for every $n).$ And if $l_{\ast}>l(\delta)$ then $2^{l_{\ast
}\varepsilon(\delta)}\geq l_{\ast}.$ By the very definition of $l_{\ast}$ we
have
\[
\frac{2^{2h(l_{\ast}-1)}}{(l_{\ast}-1)^{2}}<\theta(n_{\ast})
\]
so that
\[
2^{2hl_{\ast}}\leq2^{2h}(l_{\ast}-1)^{2}\theta(n_{\ast})\leq2^{2h}%
\times2^{2l_{\ast}\varepsilon(\delta)}\theta(n_{\ast})
\]
and, since $l(\delta)\geq1,$ this gives (\ref{reg14}).

\textbf{Step 2}. We define%
\begin{align*}
\nu_{l}  &  =0\quad if\quad l<l_{\ast}\\
&  =\mu_{n(l)}\quad if\quad l\geq l_{\ast}%
\end{align*}
and we estimate $\pi_{k,q,h,p}(\mu,(\nu_{l})_{l}).$ First, by (\ref{reg9}) and
(\ref{reg13})
\[
\sum_{l=l_{\ast}}^{\infty}\frac{1}{2^{2hl}}\left\Vert f_{n(l)}\right\Vert
_{q+2h,2h,p}\leq\sum_{l=l_{\ast}}^{\infty}\frac{1}{2^{2hl}}\theta
(n(l))\leq\Theta\sum_{l=l_{\ast}}^{\infty}\frac{1}{l^{2}}\leq\Theta.
\]
Then we write%
\[
\sum_{l=1}^{\infty}2^{(q+k+d/p_{\ast})l}d_{k}(\mu,\nu_{l})=S_{1}+S_{2}%
\]
with%
\[
S_{1}=\sum_{l=1}^{l_{\ast}-1}2^{(q+k+d/p_{\ast})l}d_{k}(\mu,0),\quad
S_{2}=\sum_{l=l_{\ast}}^{\infty}2^{(q+k+d/p_{\ast})l}d_{k}(\mu,\mu_{n(l)}).
\]
Since $d_{k}(\mu,0)\leq d_{0}(\mu,0)\leq\left\vert \mu\right\vert
({\mathbb{R}}^{d})$ we use (\ref{reg14}) and we obtain
\begin{align*}
S_{1}  &  \leq\left\vert \mu\right\vert ({\mathbb{R}}^{d})\times
2^{(q+k+d/p_{\ast})l_{\ast}}=\left\vert \mu\right\vert ({\mathbb{R}}%
^{d})\times(2^{2(h-\varepsilon(\delta))l_{\ast}})^{(q+k+d/p_{\ast
})/2(h-\varepsilon(\delta))}\\
&  \leq\left\vert \mu\right\vert ({\mathbb{R}}^{d})\times(2^{2hl(\delta
)}\theta(n_{\ast}))^{\rho_{h}(q)(1+\delta)}=A(\delta)\theta(n_{\ast}%
)^{\rho_{h}(q)(1+\delta)}.
\end{align*}
If $l\geq l_{\ast}$ then $n(l)\geq n(l_{\ast})\geq n_{\ast}$ so that, using
(\ref{reg11}) first and (\ref{reg13}) then, we obtain
\[
d_{k}(\mu,\mu_{n(l)})\leq\frac{C_{h,n_{\ast}}(\eta)}{\theta^{\rho_{h}(q)+\eta
}(n(l))}\leq C_{h,n_{\ast}}(\eta)\Big(\frac{l^{2}}{2^{2hl}}\Big)^{\rho
_{h}(q)+\eta}=\frac{C_{h,n_{\ast}}(\eta)}{2^{(q+k+d/p_{\ast})l}}\times
\frac{l^{2(\rho_{h}(q)+\eta)}}{2^{2h\eta l}}.
\]
We conclude that%
\[
S_{2}\leq C_{h,n_{\ast}}(\eta)\sum_{l=l_{\ast}}^{\infty}\frac{l^{2(\rho
_{h}(q)+\eta)}}{2^{2\eta hl}}\leq C_{h,n_{\ast}}(\eta)\times B(\eta).
\]
$\square$

We give now a consequence of the above result which is more readable.

\begin{proposition}
\label{Crit}Let $q,k,d\in\mathbb{N}$ and $p>1$ be fixed. We consider a family
of measures $\mu_{\varepsilon}(dx)=f_{\varepsilon}(x)dx,\varepsilon>0$ with
$f_{\varepsilon}\in C^{\infty}(\mathbb{R}^{d})$\ and a finite measure $\mu
$\ on $\mathbb{R}^{d}$\ which verify the following hypothesis. There exists
$\varepsilon_{\ast}>0,\beta>0,a\geq0,b\geq0,C_{0}\geq1$ and $Q_{h}(q,p)\geq1$
such that for every $h\in\mathbb{N}_{\ast}$ and every $0<\varepsilon
<\varepsilon_{\ast}$%
\begin{align}
i)\quad d_{k}(\mu_{\varepsilon},\mu)  &  \leq C_{0}\varepsilon^{\beta
},\label{reg15bb}\\
ii)\quad\left\Vert f_{\varepsilon}\right\Vert _{2h+q,2h,p}  &  \leq
Q_{h}(q,p)\varepsilon^{-b(2h+q+a)},\label{reg15c}\\
iii)\quad r  &  :=\beta-b(k+q+d/p_{\ast})>0. \label{reg15a}%
\end{align}
We denote%
\begin{equation}
h_{\ast}=\frac{1}{\varepsilon_{\ast}}\vee\frac{b(q+a)(k+q+d/p_{\ast})}{r}%
\vee\frac{q+a}{2}. \label{reg16}%
\end{equation}
Then, $\mu(dx)=f(x)dx$ with $f\in W^{q,p}(\mathbb{R}^{d}).$ Moreover, for
every $\delta>0,$ there exists a constant $C\geq1,$ depending on
$q,k,d,p,\delta,\beta,r$ and $a,b$ only (but which does not depend neither on
$h$ nor on $C_{0})$, such that for every $h\geq h_{\ast}$ one has
\begin{equation}
\left\Vert f\right\Vert _{q,p}\leq C\times C_{0}\times\left(  h^{2b}%
Q_{h}^{1/2h}(q,p)\right)  ^{(1+\delta)(k+q+d/p_{\ast})} \label{reg17}%
\end{equation}

\end{proposition}

\textbf{Proof}. All over this proof $C$ designs a constant which depends\ on
$q,k,d,p,\delta,\beta,r$ and $a,b$ only (we stress that in particular it may
depend on $C_{\ast}$ from (\ref{reg4})). We will use Lemma \ref{REG}. We take
\begin{equation}
\eta=\frac{r}{2b(2h+q+a)}\wedge(\delta\rho_{h}(q)) \label{reg17'}%
\end{equation}
with $\rho_{h}(q)$ given in (\ref{reg5}). For $h\geq h_{\ast}$ one has
$\rho_{h}(q)b(q+a)\leq\frac{r}{2},$ and, by definition, $r=\beta-2h\rho
_{h}(q)b.$ Using also (\ref{reg17'}) we obtain
\begin{align*}
\beta-(\rho_{h}(q)+\eta)b(2h+q+a)  &  =(\beta-2h\rho_{h}(q)b)-\rho
_{h}(q)b(q+a)-\eta b(2h+q+a)\\
&  \geq r-\frac{r}{2}-\frac{r}{2}=0.
\end{align*}
It follows that for every $\varepsilon\leq\varepsilon_{\ast}$ we have%
\begin{equation}
d_{k}(\mu_{\varepsilon},\mu)\left\Vert f_{\varepsilon}\right\Vert
_{2h+q,2h,p}^{\rho_{h}(q)+\eta}\leq C_{0}Q_{h}^{\rho_{h}(q)+\eta
}(q,p)\varepsilon^{\beta-(\rho_{h}(q)+\eta)b(2h+q+a)}\leq C_{0}Q_{h}^{\rho
_{h}(q)(1+\delta)}(q,p). \label{reg18}%
\end{equation}

We take now $\varepsilon_{n}=\frac{1}{n}$ and $n_{\ast}=h$ and we define
\begin{align*}
g_{n}  &  =0\quad if\quad n<n_{\ast}\\
&  =f_{\varepsilon_{n}}\quad if\quad n\geq n_{\ast}.
\end{align*}
We will use Lemma \ref{REG} for $\nu_{n}(dx)=g_{n}(x)dx$ so we have to
identify the quantities defined there. We define%
\begin{align*}
\theta(n)  &  =Q_{h}(q,p)n^{b(2h+q+a)}\quad if\quad n\geq n_{\ast},\\
\theta(n)  &  =\theta(n_{\ast})\quad if\quad n\leq n_{\ast}.
\end{align*}
By (\ref{reg15c}) \ we have $\left\Vert g_{n}\right\Vert _{2h+q,2h,p}%
\leq\theta(n)$ and moreover, for $n\geq n_{\ast}=h,$ we have%
\[
\frac{\theta(n+1)}{\theta(n)}=(1+\frac{1}{n})^{n\times\frac{b(2h+q+a)}{n}}\leq
e^{2b+\frac{b(q+a)}{h}}\leq e^{3b}.
\]
We conclude that $\Theta\leq e^{3b}.$ We estimate now $B(\eta)$ defined in
(\ref{reg12''}). Noticed first that
\begin{align*}
\frac{1}{\eta h}  &  =\frac{2b(2h+q+a)}{rh}\vee\frac{2}{\delta(q+k+d/p_{\ast
})}\\
&  \leq\frac{6b}{r}\vee\frac{2}{\delta(q+k+d/p_{\ast})}=:C_{1}%
\end{align*}
so $\eta h\geq1/C_{1}.$ Then%
\[
B(\eta)=\sum_{l=1}^{\infty}\frac{l^{2(q+k+d/p_{\ast}+\eta)}}{2^{2h\eta l}}%
\leq\sum_{l=1}^{\infty}\frac{l^{2(q+k+d/p_{\ast}+\eta)}}{2^{2l/C_{1}}}\leq C.
\]
Moreover, since $h\geq\frac{1}{2}(q+a)$ it follows that
\[
\rho_{h}(q)(1+\delta)b(2h+q+a)\leq2(1+\delta)b(k+q+d/p_{\ast})
\]
and consequently (recall that $n_{\ast}=h)$%
\begin{align*}
\theta(n_{\ast})^{\rho_{h}(q)(1+\delta)}  &  =Q_{h}^{\rho_{h}(q)(1+\delta
)}(q,p)n_{\ast}^{\rho_{h}(q)(1+\delta)b(2h+q+a)}\\
&  \leq Q_{h}^{\rho_{h}(q)(1+\delta)}(q,p)h^{2(1+\delta)b(k+q+d/p_{\ast})}.
\end{align*}
Finally we notice that, by (\ref{reg18}), the constant\ $C_{h,n_{\ast}}(\eta)$
defined in (\ref{reg11}) verifies%
\[
C_{h,n_{\ast}}(\eta)\leq C_{0}Q_{h}^{\rho_{h}(q)(1+\delta)}(q,p).
\]
As for $A(\delta)$ defined in (\ref{reg12'}), this is already a constant $C$
(which does not depend on $h$ and on $C_{0}).$ Now we use (\ref{reg12}) and we
obtain%
\[
\left\Vert f\right\Vert _{q,p}\leq C(1+Q_{h}^{\rho_{h}(q)(1+\delta
)}(q,p)h^{2b(1+\delta)(k+q+d/p_{\ast})}+C_{0}Q_{h}^{\rho_{h}(q)(1+\delta
)}(q,p))
\]
which gives (\ref{reg17}). $\square$

\subsection{Proof of Theorem \ref{Criterion}}

The aim of this section is to prove Theorem \ref{Criterion} so we consider the
framework given there: we have a family of random variables $F_{\varepsilon}$
and $G_{\varepsilon},\varepsilon>0$ such that the integration by parts formula
(\ref{ip1}) and the estimate (\ref{ip2}) hold; we also have a random variable
$F$ such that the estimate (\ref{ip3}) holds. We define the measures $\mu$ and
$\mu_{\varepsilon}$ by%
\begin{equation}
\int\phi d\mu=E(\phi(F))\quad and\quad\int\phi d\mu_{\varepsilon}%
=E(\phi(F_{\varepsilon})G_{\varepsilon}). \label{IP4}%
\end{equation}
As a consequence of (\ref{ip1}) and of (\ref{ip2}) we have $\mu_{\varepsilon
}(dx)=f_{\varepsilon}(x)dx$ with $f_{\varepsilon}\in C^{\infty}(\mathbb{R}%
^{d}).$

We also consider a function $\Phi\in\mathcal{C}$ (so in particular $\Phi$
verifies (\ref{ip4})). All these hypothesis are in force in this section.

Moreover, for $v\in\mathbb{R}^{d}$ we construct the "exterior rectangle"
$A_{v}$ in the following way. For $y\in\mathbb{R}$ we denote $I_{y}%
=(y,\infty)$ if $y\geq0$ and $I_{y}=(-\infty,y)$ if $y<0.$ And for
$v=(v_{1},...,v_{d})$ we define
\begin{equation}
A_{v}=\prod_{i=1}^{d}I_{v_{i}}. \label{IP2a}%
\end{equation}

We will first prove the following two Lemmas.

\begin{lemma}
\label{NORM} For every $q,h\in\mathbb{N}$ and $p>1$ there exists some
constants $C$ and $\theta$ (depending on $q,h,d$ and $p$ only) and
$\varepsilon_{\ast}>0$ such that, for every $\varepsilon\in(0,\varepsilon
_{\ast})$%

\begin{equation}
\left\Vert \Phi f_{\varepsilon}\right\Vert _{q,h,p}\leq C\times C_{\theta
}(\Phi)\times\widehat{H}_{q+d,p_{\ast}}\times\varepsilon^{-b(q+d+a)}
\label{IP5}%
\end{equation}
with $C_{\theta}(\Phi)$ given in (\ref{ip5}) and $\widehat{H}_{q+d,p_{\ast}}%
$\ given in (\ref{ip2}).
\end{lemma}

\textbf{Proof}. We denote%
\[
I_{q,h,p}(\Phi)(x)=\sup_{\left\vert \beta\right\vert \leq q}\int_{R^{d}%
}(1+\left\vert v\right\vert )^{h}\left\vert \partial^{\beta}\Phi(v)\right\vert
^{p}1_{A_{v}}(x)dv.
\]

By (\ref{ip1}) (we use a formal computation which may be done rigorous by
regularization by convolution)
\begin{align*}
\partial^{\alpha}(\Phi f_{\varepsilon})(v)  &  =\sum_{(\beta,\gamma)=\alpha
}\partial^{\beta}\Phi(v)\partial^{\gamma}f_{\varepsilon}(v)=\sum
_{(\beta,\gamma)=\alpha}\partial^{\beta}\Phi(v)E(\partial^{\gamma}\delta
_{0}(F_{\varepsilon}-v)G_{\varepsilon})\\
&  =\sum_{(\beta,\gamma)=\alpha}\partial^{\beta}\Phi(v)E(1_{A_{v}%
}(F_{\varepsilon})H_{(\gamma,1,...,d),\varepsilon}).
\end{align*}
Using H\"{o}lder's inequality first and then (\ref{ip2})%
\begin{align*}
\left\vert \partial^{\alpha}(\Phi f_{\varepsilon})(v)\right\vert  &  \leq
\sum_{(\beta,\gamma)=\alpha}\left\vert \partial^{\beta}\Phi(v)\right\vert
P^{1/p}(F_{\varepsilon}\in A_{v})\left\Vert H_{(\gamma,1,...,d),\varepsilon
}\right\Vert _{p_{\ast}}\\
&  \leq\sum_{\left\vert \beta\right\vert \leq q}\left\vert \partial^{\beta
}\Phi(v)\right\vert P^{1/p}(F_{\varepsilon}\in A_{v})\times\widehat
{H}_{q+d,p_{\ast}}\varepsilon^{-b(q+d+a)}.
\end{align*}
This gives%
\begin{align*}
\left\Vert \Phi f_{\varepsilon}\right\Vert _{q,h,p}  &  =\sum_{\left\vert
\alpha\right\vert \leq q}(\int(1+\left\vert v\right\vert )^{h}\left\vert
\partial^{\alpha}(\Phi f_{\varepsilon})(v)\right\vert ^{p}dv)^{1/p}\\
&  \leq C\widehat{H}_{q+d,p_{\ast}}\varepsilon^{-b(q+d+a)}\sum_{\left\vert
\beta\right\vert \leq q}[\int(1+\left\vert v\right\vert )^{h}\left\vert
\partial^{\beta}\Phi(v)\right\vert ^{p}E(1_{F_{\varepsilon}\in A_{v}%
})dv]^{1/p}\\
&  =C\widehat{H}_{q+d,p_{\ast}}\varepsilon^{-b(q+d+a)}[E(I_{q,h,p}%
(\Phi)(F_{\varepsilon}))]^{1/p}.
\end{align*}
By using (\ref{ip4}) we may find some constants $c_{1}$ and $c_{2}$ such that
for every $\gamma$ with $\left\vert \gamma\right\vert \leq q$%
\[
\left\vert \partial^{\gamma}\Phi(v)\right\vert ^{p}\leq\frac{c_{1}%
}{(1+\left\vert v\right\vert )^{h+d+1}}\left\vert \Phi(v)\right\vert ^{c_{2}p}%
\]
so we get
\[
E(I_{q,h,p}(\Phi)(F_{\varepsilon}))\leq c_{1}E(\int_{R^{d}}(1+\left\vert
v\right\vert )^{-(d+1)}\left\vert \Phi(v)\right\vert ^{c_{2}p}1_{A_{v}%
}(F_{\varepsilon})dv).
\]
If $x\in A_{v}$ then $\left\vert v\right\vert \leq\left\vert x\right\vert $
and since $\Phi\in\mathcal{C}$ this implies $\Phi(v)\leq C\Phi(x).$ So the
above term is upper bounded by%
\[
CE(\int_{R^{d}}(1+\left\vert v\right\vert )^{-(d+1)}\left\vert \Phi
(F_{\varepsilon})\right\vert ^{c_{2}p}1_{A_{v}}(F_{\varepsilon})dv)\leq
CE(\left\vert \Phi(F_{\varepsilon})\right\vert ^{c_{2}p})\leq C\times
C_{\theta}(\Phi)
\]
with $\theta=c_{2}p.\ \square$

\begin{lemma}
\label{APROXIMATION}We recall that by (\ref{ip3})%
\begin{equation}
\left\Vert 1-G_{\varepsilon}\right\Vert _{2}+\left\Vert F-F_{\varepsilon
}\right\Vert _{1}\leq C_{\ast}\varepsilon^{\beta}. \label{IP8}%
\end{equation}
Then, for every $\delta>0$ there exists $\theta(\delta)\geq1$ and $C\geq
1$\ such that
\begin{equation}
d_{1}(\Phi\mu_{\varepsilon},\Phi\mu)\leq C(C_{\ast}+C_{\theta(\delta)}%
(\Phi)))\varepsilon^{\beta(1-\delta)} \label{IP9'}%
\end{equation}
with $C_{\theta}(\Phi)$ defined in (\ref{ip5}).
\end{lemma}

\textbf{Proof. }Let $\phi$ with $\left\Vert \phi\right\Vert _{1,\infty}\leq1.$
We estimate first%
\begin{equation}
\left\vert E((\phi\Phi)(F_{\varepsilon})(1-G_{\varepsilon}))\right\vert
\leq\left\Vert \phi\right\Vert _{\infty}\left\Vert \Phi(F_{\varepsilon
})\right\Vert _{2}\left\Vert 1-G_{\varepsilon}\right\Vert _{2}\leq C_{2}%
(\Phi)\times C_{\ast}\varepsilon^{\beta}. \label{IP14'}%
\end{equation}
Then we write
\begin{equation}
\left\vert E((\phi\Phi)(F_{\varepsilon})-(\phi\Phi)(F))\right\vert \leq
E\int_{0}^{1}\left\vert \nabla(\phi\Phi)(\lambda F+(1-\lambda)F_{\varepsilon
})(F-F_{\varepsilon})\right\vert d\lambda. \label{IP14}%
\end{equation}
Using (\ref{ip4}) and the fact that $\Phi$ is a convex function
\begin{align*}
\left\vert \nabla(\phi\Phi)(\lambda F+(1-\lambda)F_{\varepsilon})\right\vert
&  \leq c_{1}\left\Vert \phi\right\Vert _{1,\infty}\left\vert \Phi(\lambda
F+(1-\lambda)F_{\varepsilon})\right\vert ^{c_{2}}\\
&  \leq C(\lambda\left\vert \Phi(F)\right\vert ^{c_{2}}+(1-\lambda)\left\vert
\Phi(F_{\varepsilon})\right\vert ^{c_{2}}).
\end{align*}
It follows that the last term in (\ref{IP14}) is upper bounded by%
\[
C(E(\left\vert \Phi(F)\right\vert ^{c_{1}}\left\vert F-F_{\varepsilon
}\right\vert )+E(\left\vert \Phi(F_{\varepsilon})\right\vert ^{c_{1}%
}\left\vert F-F_{\varepsilon}\right\vert )).
\]
We take $K>0$ and we write
\[
E(\left\vert \Phi(F_{\varepsilon})\right\vert ^{c_{1}}\left\vert
F-F_{\varepsilon}\right\vert )=I_{K}(F_{\varepsilon})+J_{K}(F_{\varepsilon})
\]
with%
\[
I_{K}(F_{\varepsilon})=E(\left\vert \Phi(F_{\varepsilon})\right\vert ^{c_{1}%
}\left\vert F-F_{\varepsilon}\right\vert 1_{\{\left\vert \Phi(F_{\varepsilon
})\right\vert ^{c_{1}}\leq K\}})\leq KC_{\ast}\varepsilon^{\beta}%
\]
and%
\begin{align*}
J_{K}(F_{\varepsilon})  &  =E(\left\vert \Phi(F_{\varepsilon})\right\vert
^{c_{1}}\left\vert F-F_{\varepsilon}\right\vert 1_{\{\left\vert \Phi
(F_{\varepsilon})\right\vert ^{c_{1}}>K\}})\\
&  \leq(\left\Vert F\right\Vert _{2}+\left\Vert F_{\varepsilon}\right\Vert
_{2})(E(\left\vert \Phi(F_{\varepsilon})\right\vert ^{2c_{1}}1_{\{\left\vert
\Phi(F_{\varepsilon})\right\vert ^{c_{1}}>K\}}))^{1/2}\\
&  \leq\frac{(\left\Vert F\right\Vert _{2}+\left\Vert F_{\varepsilon
}\right\Vert _{2})}{K^{\theta}}(E(\left\vert \Phi(F_{\varepsilon})\right\vert
^{2c_{1}(1+\theta)})^{1/2}.
\end{align*}
By (\ref{ip4}) $\left\Vert F\right\Vert _{2}+\left\Vert F_{\varepsilon
}\right\Vert _{2}\leq C\times C_{c_{2}}^{1/2}(\Phi)$ so finally%
\[
J_{K}(F_{\varepsilon})\leq\frac{C}{K^{\theta}}C_{c_{2}\vee2c_{1}(1+\theta
)}(\Phi).
\]
These estimates hold for every $\theta\geq1$ and $K>0.$ In order to optimize
we take $K=\varepsilon^{-\beta/(1+\theta)}$ and we obtain%
\[
I_{K}(F_{\varepsilon})+J_{K}(F_{\varepsilon})\leq C\times C_{c_{2}\vee
2c_{1}(1+\theta)}(\Phi)\times\varepsilon^{\beta\times\frac{\theta}{1+\theta}%
}.
\]
A similar inequality holds with $F$ instead of $F_{\varepsilon},$ so, using
(\ref{IP14'}) as well we obtain
\[
\left\vert E((\phi\Phi)(F_{\varepsilon})-(\phi\Phi)(F))\right\vert \leq
C(C_{\ast}+C_{c_{2}\vee2c_{1}(1+\theta)}(\Phi))\times\varepsilon^{\beta
\times\frac{\theta}{1+\theta}}%
\]
Then, taking $\theta=(1-\delta)/\delta$ we get (\ref{IP9'}). $\square$

We are now ready to give:

\textbf{Proof of Theorem} \ref{Criterion}. We will use the Proposition
\ref{Crit} with $k=1$ for the measures $(\Phi\mu)(dx)$ and $(\Phi
\mu_{\varepsilon})(dx)=\Phi(x)f_{\varepsilon}(x)dx$ with $\mu$ and
$\mu_{\varepsilon}$ given in (\ref{IP4}). By (\ref{ip6}) we may find (and fix)
$\delta>0$ such that%
\begin{equation}
b-\delta>b(1+q+\frac{d}{p_{\ast}}). \label{b}%
\end{equation}
By (\ref{IP9'}) we have%
\[
d_{1}(\Phi\mu_{\varepsilon},\Phi\mu)\leq C(C_{0}+C_{\theta(\delta)}%
(\Phi)))\varepsilon^{\beta(1-\delta)}%
\]
so (\ref{reg15bb}) holds with the constant $C_{0}=C(C_{\ast}+C_{\theta
(\delta)}(\Phi)))$. By (\ref{IP5})
\[
\left\Vert \Phi f_{\varepsilon}\right\Vert _{2h+q,2h,p}\leq C\times C_{\theta
}(\Phi)\times\widehat{H}_{2h+q+d,p_{\ast}}\varepsilon^{-b(2h+q+d+a)}.
\]
So (\ref{reg15c}) holds with $Q_{h}(q,p)=CC_{\theta}(\Phi)\widehat
{H}_{2h+q+d,p_{\ast}}.$ Notice that $a$ from (\ref{reg15c}) is replaced by
$a^{\prime}=a+d.$ This changes nothing except the value of $h_{\ast}$ (see
(\ref{reg16})) which anyway, is not explicit in our statement. As a
consequence of (\ref{b}) hypothesis (\ref{reg15a}) is verified so we are able
to use Proposition \ref{Crit} in order to obtain $(\Phi\mu)(dx)=g(x)dx$ with
$g\in W^{q,p}(\mathbb{R}^{d}).$ It follows that $\mu(dx)=f(x)dx$ with
$f=g/\Phi$ and by (\ref{reg17}) (the value of $\theta$ changes from a line to
another and we use the inequality $C_{\theta}^{\theta^{\prime}}(\Phi)\leq
C_{\theta\times\theta^{\prime}}(\Phi)$)
\begin{align*}
\left\Vert \Phi f\right\Vert _{q,p}  &  \leq C\times(C_{\ast}+C_{\theta}%
(\Phi))\times\left(  h^{2b}Q_{h}^{1/2h}(q,p)\right)  ^{(1+\delta
)(k+q+d/p_{\ast})}\\
&  \leq C\times(C_{\ast}+C_{\theta}(\Phi))\times\left(  h^{2b}C_{\theta
}^{1/2h}(\Phi)\widehat{H}_{2h+q+d,p_{\ast}}^{1/2h}\right)  ^{(1+\delta
)(k+q+d/p_{\ast})}\\
&  \leq\Gamma_{\Phi,\theta}(q,h,p).
\end{align*}
So (\ref{ip7}) is proved.

We prove now the point \textbf{C}. If the above inequality holds with $q=1$
and $p>d$, then, by Morrey's lemma $\Phi f$ is $\chi-$ H\"{o}lder continuous
with $\chi=1-\frac{d}{p}$ and
\[
\left\Vert \Phi f\right\Vert _{\infty}\leq\left\Vert \Phi f\right\Vert
_{C^{0,\chi}}\leq C\left\Vert \Phi f\right\Vert _{1,p}\leq C\Gamma_{\Phi
}(1,h,p).
\]
So we obtain (\ref{ip9}). Let us prove (\ref{ip11}). Recall the definition of
$\widehat{\nabla}\Phi(x)$ defined in (\ref{ip10}). We write%
\[
(\Phi f)(y)-(\Phi f)(x)=\Phi(x)(f(y)-f(x))+(\Phi(y)-\Phi(x))f(x)
\]
so that, if $\left\vert x-y\right\vert \leq1,$ then
\begin{align*}
\left\vert \Phi(x)(f(y)-f(x))\right\vert  &  \leq\left\Vert \Phi f\right\Vert
_{C^{0,\chi}}\left\vert x-y\right\vert ^{\chi}+\widehat{\nabla}\Phi
(x)\left\vert f(x)\right\vert \left\vert x-y\right\vert \\
&  \leq C\left\Vert \Phi f\right\Vert _{1,p}\left\vert x-y\right\vert ^{\chi
}+\frac{\widehat{\nabla}\Phi(x)}{\Phi(x)}\left\Vert \Phi f\right\Vert
_{1,p}\left\vert x-y\right\vert \\
&  \leq C\Gamma_{\Phi}(1,h,p)\left\vert x-y\right\vert ^{\chi}+\widehat
{\nabla}\Phi(x)\times\frac{\Gamma_{\Phi}(1,h,p)}{\Phi(x)}\left\vert
x-y\right\vert
\end{align*}
where, in order to obtain the second inequality we have used (\ref{ip9}). So
(\ref{ip11}) is proved. $\square$

\bigskip

\section{\bigskip Appendix 2: The integration by parts formula.}

The aim of this section is to give a hint to the proof of (\ref{bo9a}) and of
the estimates (\ref{bo9b}) and (\ref{bo9}). The integration by parts formula
(\ref{bo9a}) has been established in \cite{[BF]} by using a version of
Malliavin calculus for jump processes introduced in \cite{[BCl1]}. All this
machinery is quit heavy and we are not able to give here a detailed technical
view (we refer to \cite{[BF]} for a complete presentation). We just try to
give an overlook which permits to the reader to understand which are the main
objects and arguments involved in the proof of these results. And also to
precise the dependence with respect to $t$ of the constants in (\ref{bo9b})
and (\ref{bo9}).

\subsection{Real chock and fictive chock representation}

A first step is to give some appropriate alternative representations for
$V_{t}^{\varepsilon,\zeta},$ solution of the equation (see (\ref{bo7'})):%

\begin{equation}
V_{t}^{\varepsilon,\zeta}=V_{0}+\int_{0}^{t}\int_{E\times\mathbb{R}_{+}%
}A(\theta)(V_{s-}^{\varepsilon,\zeta}-v)1_{\{u\leq\varphi_{\varepsilon
}^{\gamma}(\left\vert V_{s-}^{\varepsilon,\zeta}-v\right\vert )\}}I_{\zeta
}(\theta)N(ds,d\theta,dv,du). \label{a1}%
\end{equation}
We recall $E=[-\frac{\pi}{2},\frac{\pi}{2}]\times\mathbb{R}^{2}$ and
$N(dt,d\theta,dv,du)$ is a Poisson point measure on $E\times\mathbb{R}_{+}$
with intensity measure $dt\times b(\theta)d\theta\times df_{t}(dv)\times du$
where $f_{t}(dv)$ is the solution (which exists and is unique) of the equation
(\ref{bo1}).

\textbf{Step 1}. In a first stage we use some change of variable in order to
write the above equation in an alternative form which is appropriate for our
calculus (see Section 3 in \cite{[BF]} for details; the motivation of this new
representation is just technical). Using the Skorohod representation theorem
we may find a measurable function $v_{t}:[0,1]\rightarrow\mathbb{R}^{2}$ such
that for every $\psi:\mathbb{R}^{2}\rightarrow\mathbb{R}_{+}$%
\[
\int_{0}^{1}\psi(v_{t}(\rho))d\rho=\int_{\mathbb{R}^{2}}\psi(v)f_{t}(dv).
\]
This allows to replace the measure $f_{t}(dv)$ on $\mathbb{R}^{2}$ by $d\rho$
on $[0,1].$

Moreover, for $x\in(0,\frac{\pi}{2}],$ let $G(x)=\int_{x}^{\pi/2}%
b(\theta)d\theta$ and let $g:(0,\infty)\rightarrow(0,\frac{\pi}{2}]$ be the
inverse of $G,$ that is $G(g(z))=z$ (since $b(\theta)\simeq\left\vert
\theta\right\vert ^{-(1+\nu)}$ by assumption, it follows that $G(x)\simeq
\nu^{-1}(x^{-\nu}-(\pi/2)^{-\nu})$ and $g(z)\simeq(1+z)^{-\nu}).$ For
$z<0$\ we define $g(z)=-g(-z)$. With this construction we will have%
\[
\int_{-\pi/2}^{\pi/2}\psi(\theta)b(\theta)d\theta=\int_{\mathbb{R}_{\ast}%
}^{\infty}\psi(g(z))dz
\]
and this allows to replace the measure $b(\theta)d\theta$ on $(-\frac{\pi}%
{2},\frac{\pi}{2})$ by $dz$ on $\mathbb{R}_{\ast}:=\mathbb{R}\diagdown\{0\}.$

Finally we consider a function $\mathbf{I}_{\zeta}:\mathbb{R}\rightarrow
\lbrack0,1]$ which is smooth, with all derivatives bounded and such that
$\mathbf{I}_{\zeta}(z)=1$ for $\left\vert z\right\vert \leq G(\zeta)$ and
$\mathbf{I}_{\zeta}(z)=0$ for $\left\vert z\right\vert \geq G(\zeta)+1.$ And
we choose the function $I_{\zeta}$ in equation (\ref{a1}) in such a way that%
\begin{equation}
\mathbf{I}_{\zeta}(z)=I_{\zeta}(g(z)). \label{a2}%
\end{equation}
Then we may write the equation (\ref{a1}) as
\begin{equation}
V_{t}^{\varepsilon,\zeta}=V_{0}+\int_{0}^{t}\int_{0}^{1}\int_{-G(\zeta
)-1}^{G(\zeta)+1}\int_{0}^{2\Gamma_{\varepsilon}^{\gamma}}A(g(z))(V_{s-}%
^{\varepsilon,\zeta}-v_{s}(\rho))1_{\{u\leq\varphi_{\varepsilon}^{\gamma
}(\left\vert V_{s-}^{\varepsilon,\zeta}-v_{s}(\rho)\right\vert )\}}%
\mathbf{I}_{\zeta}(z)M(ds,d\rho,dz,du) \label{a3}%
\end{equation}
with $M$ a Poisson point measure on $[0,T]\times\lbrack0,1]\times
\mathbb{R}_{\ast}\times(0,\infty)$ with intensity measure $m(ds,d\rho
,dz,du)=dsd\rho dzdu.$ Notice that we may take $u\leq2\Gamma_{\varepsilon
}^{\gamma}$ because we know that $\varphi_{\varepsilon}(v)\leq\Gamma
_{\varepsilon}$ (we use $2\Gamma_{\varepsilon}$ instead of $\Gamma
_{\varepsilon}$ just for technical reasons)$.$

\textbf{Step 2}. Since $m$ is a finite measure we may represent the above
equation by using a compound Poisson process as follows. We consider a
standard Poisson process $J_{t}^{\varepsilon}=\sum_{k=1}^{\infty}%
1_{\{T_{k}\leq t\}}$ of parameter $\lambda=4(G(\zeta)+1)\Gamma_{\varepsilon
}^{\gamma}$ and a sequence $(\overline{R}_{k},\overline{Z}_{k},\overline
{U}_{k}),k\in\mathbb{N}$ of independent random variables, uniformly
distributed on $[0,1]\times\lbrack-G(\zeta)-1,G(\zeta)+1]\times\lbrack
0,2\Gamma_{\varepsilon}^{\gamma}],$ which are independent of $J$ (all these
objects depend on $\varepsilon$ and $\zeta,$ but, as they are fixed, we do not
mention it in the notation). Now the equation (\ref{a3}) may be written as
\begin{equation}
V_{t}^{\varepsilon,\zeta}=V_{0}+\sum_{T_{k}\leq t}A(g(\overline{Z}%
_{k}))(V_{T_{k-1}}^{\varepsilon,\zeta}-v_{T_{k}}(\overline{R}_{k}%
))1_{\{\overline{U}_{k}\leq\varphi_{\varepsilon}^{\gamma}(\left\vert
V_{T_{k-1}}^{\varepsilon,\zeta}-v_{T_{k}}(\overline{R}_{k})\right\vert
)\}}\mathbf{I}_{\zeta}(\overline{Z}_{k}). \label{a4}%
\end{equation}
This equation is known as "the fictive chock" representation (see for example
\cite{PSL}).

\textbf{Step 3}. The idea of the approach by means of the Malliavin calculus
is to look to $V_{t}^{\varepsilon,\zeta}$ as a functional $f(\overline{Z}%
_{1},...,\overline{Z}_{J_{t}})$ of $\overline{Z}_{k},k=1,...,J_{t}%
^{\varepsilon}$ and to use an elementary integration by parts formula based on
the (uniform) law of $\overline{Z}_{k}.$ But this is not possible directly
because of the indicator function of the set $\{\overline{U}_{k}\leq
\varphi_{\varepsilon}^{\gamma}(\left\vert V_{T_{k-1}}^{\varepsilon,\zeta
}-v_{T_{k}}(\overline{R}_{k})\right\vert )\}$ which appears in the equation
:\ as a consequence $\overline{Z}_{k}\rightarrow f(\overline{Z}_{1}%
,...,\overline{Z}_{J_{t}^{\varepsilon}})=V_{T_{k}}^{\varepsilon,\zeta}$ is not
differentiable. In order to avoid this difficulty we introduce the so called
"real chock representation" that we present now. We consider the equation%
\begin{equation}
V_{t}^{\varepsilon,\zeta}=V_{0}+\sum_{T_{k}\leq t}A(g(Z_{k}))(V_{T_{k-1}%
}^{\varepsilon,\zeta}-v_{T_{k}}(R_{k}))\mathbf{I}_{\zeta}(Z_{k}). \label{a5}%
\end{equation}
This is the same equation as (\ref{a4}) but the indicator function disappears.
But now the law of $(R_{k},Z_{k})$ is no more the uniform law (as above). We
define this law as follows: we assume that conditionally to $T_{k}=t$ and
$V_{T_{k}}^{\varepsilon,\zeta}=w$ the law of $(R_{k},Z_{k})$ is given by%
\[
P_{t,w}(R_{k}\in d\rho,Z_{k}\in dz)=q_{t,w}(\rho,z)dzd\rho
\]
with%
\[
q_{t,w}(\rho,z)=\frac{1}{\lambda}\varphi_{\varepsilon}^{\gamma}(\left\vert
w-v(\rho)\right\vert )1_{\{\left\vert z\right\vert \leq G(\zeta)+1\}}%
+g_{t,w}(z).
\]
Here $g_{t,w}(z)$ is an auxiliary smooth function which is null on
$\{\left\vert z\right\vert \leq G(\zeta)+1\}$ and which is chosen in such a
way that $\int\int q_{t,w}(\rho,z)d\rho dz=1$ (it plays the role of a cemetery
and does not come on in the computations). Notice that $q_{t,w}(\rho
,z)dzd\rho$ gives the precise way in which the law of $(R_{k},Z_{k})$ (and so
the law of the jump) depends on the position $V_{T_{k}-}^{\varepsilon,\zeta
}=w$ of the particle. One may check (see Section 3 in \cite{[BF]} for details)
that the law of the solution of the equation (\ref{a5})\ coincides with the
law of $V_{t}^{\varepsilon,\zeta},$ the solution of (\ref{a4}). So we may (and
do) work with $V_{t}^{\varepsilon,\zeta}$ solution of (\ref{a5}) now on. This
is the "real chock representation". Now the machinery which produces
$V_{t}^{\varepsilon,\zeta}$ as a function of $Z_{k}$ is a smooth function and
we may use a differential calculus for it. Notice that we know nothing about
the regularity of the function $\rho\rightarrow v_{t}(\rho)$ and consequently
we are not able to use $R_{k}$ - we will just use $Z_{k}.$ We also mention
that $V_{T_{k}}^{\varepsilon,\zeta}$ is a function of $T_{i},R_{i}%
,Z_{i},i=1,...,k$ so we will use the (slightly abusive notation)%
\begin{equation}
V_{T_{k}}^{\varepsilon,\zeta}=\mathcal{H}_{k}(\omega,Z_{1},...,Z_{k})
\label{a6}%
\end{equation}
where $\omega$ indicates the dependence on $T_{i},R_{i},i=1,...,k.$

\subsection{Finite dimensional Malliavin calculus.}

In this section we present the results concerning the Malliavin calculus based
on the random variables $Z_{k},k\in\mathbb{N}$ from the previous section. We
add two standard normal random variables $Z_{-1},Z_{0}$ which are independent
of $Z_{k},k=1,2,...$ as well (they correspond to the two dimensional standard
normal random variable $Z$ in introduced (\ref{bo8d})). Given $t>0$ we denote
$Z_{t}=(Z_{-1},Z_{0},Z_{1},...,Z_{J_{t}^{\varepsilon}}).$ The law of $Z_{t}$
is absolutely continuous with respect to the Lebesgue measure on
$\mathbb{R}^{J_{t}^{\varepsilon}+2}$ and has the density%
\begin{equation}
p_{t}(\omega,z)=c_{t}e^{\frac{\left\vert z_{-1}\right\vert ^{2}+\left\vert
z_{0}\right\vert ^{2}}{2}}\prod_{k=1}^{J_{t}^{\varepsilon}}q_{T_{k}%
,\mathcal{H}_{k-1}(\omega,z_{1},...,z_{k-1})}(R_{k},z_{k}). \label{a9}%
\end{equation}
The integration by parts formula which we derive in the sequel will be based
on the logarithmic derivative of this density. In order to avoid border terms
in the integration by parts formula we introduce the weights
\[
\pi_{-1}=\pi_{0}=1\quad and\quad\pi_{k}=a_{\zeta}(Z_{k})
\]
where $a_{\zeta}$ is a smooth version of $1_{(1,G(\zeta))}(z).$

We follow the strategy established in Malliavin calculus. A simple functional
$F$ is a random variable of the form $F=h(\omega,Z_{t})$ where $\omega$
designees the dependence on $T_{i},R_{i},i\in\mathbb{N}$ and $z\rightarrow
h(\omega,z)$ is a smooth function on $\mathbb{R}^{J_{t}^{\varepsilon}+2}$
which has bounded derivatives of any order. Then we define the derivatives
\[
D_{k}F=\pi_{k}\partial_{z_{k}}h(\omega,Z_{t})
\]
and for a multi-index $\kappa=(k_{1},...,k_{m})\in\{-1,0,1,...,J_{t}%
^{\varepsilon}\}^{m}$ we denote $\left\vert \kappa\right\vert =m$\ and we
define
\[
D^{\kappa}F=D_{k_{m}}...D_{k_{1}}F.
\]
For a $d$ dimensional simple functional $F=(F_{1},...,F_{d})$ the Malliavin
covariance matrix is given by%
\[
\sigma_{F}^{i,j}=\left\langle DF_{i},DF_{j}\right\rangle =\sum_{k=-1}%
^{J_{t}^{\varepsilon}}D_{k}F_{i}\times D_{k}F_{j}%
\]
and, for a one dimensional simple functional $F,$ the divergence operator is
defined by%
\begin{equation}
LF=-\sum_{k=-1}^{J_{t}^{\varepsilon}}\left[  \frac{1}{\pi_{k}}D_{k}(\pi
_{k}D_{k}F)+D_{k}F\times D_{k}\log p_{t}(\omega,Z_{t})\right]  . \label{a16}%
\end{equation}
Using elementary integration by parts on $\mathbb{R}^{J_{t}^{\varepsilon}+2}$
one obtains the following duality formula%
\[
E(\left\langle DF,DG\right\rangle )=E(FLG)=E(GLF).
\]
We will work with the norms%
\[
\left\vert F\right\vert _{m}=\left\vert F\right\vert +\sum_{1\leq\left\vert
\kappa\right\vert \leq m}\left\vert D^{\kappa}F\right\vert .
\]
The standard arguments from Malliavin calculus give the following integration
by parts formula (see Theorem 1 and 3 in \cite{[BCl1]})): let $G$ and
$F=(F_{1},...,F_{d})$ be simple functionals. We suppose that $1/\det\sigma
_{F}\in\cap_{p>1}L^{p}$. Then for every $\psi\in C_{b}^{\infty}(R^{d})$ and
every multi-index $\beta=(\beta_{1},...,\beta_{q})\in\{1,...,d\}^{q}$ one has
\begin{equation}
E(\partial^{\beta}\psi(F)G))=E(\psi(F)K_{\beta}(F,G)) \label{a7}%
\end{equation}
where $K_{\beta}(F,G)$ is a random variable which verifies%
\begin{align}
\left\vert K_{\beta}(F,G)\right\vert  &  \leq C\times\overline{K}_{\beta
}(F,G)\times\left\vert \det\sigma_{F}\right\vert ^{-(3q-1)}\quad
with\label{a8}\\
\overline{K}_{\beta}(F,G)  &  =\left\vert G\right\vert _{q}(1+\left\vert
F\right\vert _{q+1})^{q(6d+1)}\left(  1+\sum_{j=1}^{q}\sum_{k_{1}%
+...+k_{j}=q-j}\prod_{i=1}^{j}\left\vert LF\right\vert _{k_{i}}\right)  .
\label{a8'}%
\end{align}

In our approach (see Section 3) we choose $\zeta(\varepsilon)=\varepsilon
^{(1+\gamma+\alpha)/(1-\nu)}$ and we use above estimate for (see
(\ref{bo8d})):
\[
F_{t}^{\varepsilon,\zeta(\varepsilon)}=\sqrt{u_{\zeta(\varepsilon)}(t)}%
Z+V_{t}^{\varepsilon,\zeta(\varepsilon)}\ \quad with\quad u_{\zeta}%
(t)=t\zeta^{4+\nu}.
\]
We also recall that in (\ref{bo8c}) we have introduced $G_{t}^{\varepsilon
,\zeta}$ which is a smooth version of the indicator function of the set
$\{\sup_{s\leq t}\left\vert V_{s}^{\zeta,\varepsilon}\right\vert \leq
\Gamma_{\varepsilon}\}.$ In particular, for every $q,$ we have $\left\vert
G_{t}^{\varepsilon,\zeta}\right\vert _{q}=0$ on $\{\sup_{s\leq t}\left\vert
V_{s}^{\zeta,\varepsilon}\right\vert >\Gamma_{\varepsilon}\}.$ It follows
that
\begin{equation}
\overline{K}_{\beta}(F_{t}^{\varepsilon,\zeta(\varepsilon)},G_{t}%
^{\varepsilon,\zeta(\varepsilon)})=\overline{K}_{\beta}(F_{t}^{\varepsilon
,\zeta(\varepsilon)},G_{t}^{\varepsilon,\zeta(\varepsilon)})1_{\{\sup_{s\leq
t}\left\vert V_{s}^{\zeta,\varepsilon}\right\vert \leq\Gamma_{\varepsilon}\}}.
\label{a8''}%
\end{equation}

\begin{remark}
\label{D} The main difficulty in our approach comes from the estimate of
$LF_{t}^{\varepsilon,\zeta(\varepsilon)}$ which blows up as $\varepsilon
\rightarrow0.$ In order to understand this we stress that the definition of
$LF_{t}^{\varepsilon,\zeta(\varepsilon)}$\ involves%
\[
\partial_{z_{k}}\log p_{\varepsilon,t}(\omega,z)=\sum_{i=k}^{J_{t}%
^{\varepsilon}}\partial_{z_{r}}\log q_{T_{i},\mathcal{H}_{i-1}(\omega
,z_{1},...,z_{i-1})}(R_{i},z_{i}).
\]
If $q_{t,w}(\rho,z)$ does not depend on $w$ (this means that the law of the
jump does not depend on the position of the particle) then only the first term
corresponding to $i=k$ in the sum is non null. But if it does depend (and this
is our case), then all the terms are non null because of $\mathcal{H}%
_{i-1}(\omega,z_{1},...,z_{i-1})$ depends on $z_{k}$ for every
$i=k+1,...,J_{i}^{\varepsilon}$ (the perturbation of $z_{k}$ propagates in the
future)$.$ So we have $J_{t}^{\varepsilon}$\ terms in the sum and, since
$E(J_{t}^{\varepsilon})\rightarrow\infty$\ as $\varepsilon\rightarrow0,$ this
generates a blow-up and we have to give an accurate estimate of it. It
represents the main difficulty in our approach.
\end{remark}

Our aim now is to give an upper bound for the $L^{p}$ norm of $K_{\beta}%
(F_{t}^{\varepsilon,\zeta},G_{t}^{\varepsilon,\zeta})$ more precisely (see
(\ref{bo9b}))%
\begin{equation}
\left\Vert K_{\beta}(F_{t}^{\varepsilon,\zeta},G_{t}^{\varepsilon,\zeta
})\right\Vert _{p}\leq\frac{Ce^{C\Gamma_{\varepsilon}^{\gamma}}}%
{t^{4(3q-1)\times\frac{2+\nu}{\nu}}}(\varepsilon^{-q}\zeta^{-\nu q}%
+e^{-\Gamma_{\varepsilon}^{\kappa}}\zeta^{-2\nu q}) \label{a13}%
\end{equation}
with $q=\left\vert \beta\right\vert .$

\textbf{Sketch of the proof.} In Proposition 4.11 in \cite{[BF]} one proves
that for each $p>1,l\in\mathbb{N}$ and $T>0$
\begin{align*}
E(1_{\{\sup_{[0,T]}\left\vert V_{s}^{\varepsilon,\zeta}\right\vert \leq
\Gamma_{\varepsilon}\}}\sup_{[0,T]}\left\vert V_{s}^{\varepsilon,\zeta
}\right\vert _{l}^{p})  &  \leq Ce^{C\Gamma_{\varepsilon}^{\gamma}}\quad and\\
E(1_{\{\sup_{[0,T]}\left\vert V_{s}^{\varepsilon,\zeta}\right\vert \leq
\Gamma_{\varepsilon}\}}\sup_{[0,T]}\left\vert LV_{s}^{\varepsilon,\zeta
}\right\vert _{l}^{p})  &  \leq\frac{Ce^{C\Gamma_{\varepsilon}^{\gamma}}%
}{\varepsilon^{p(l+1)}\zeta^{\nu p}}%
\end{align*}
where $C$ is a constant which depends on $p,l$ and $T.$ Using the above
estimates (recall (\ref{a8''})) and H\"{o}lder's inequality one obtains%
\begin{equation}
\sup_{t\leq T}\left\Vert \overline{K}_{\beta}(F_{t}^{\varepsilon,\zeta}%
,G_{t}^{\varepsilon,\zeta})\right\Vert _{p}\leq Ce^{C\Gamma_{\varepsilon
}^{\gamma}}(\varepsilon^{-q}\zeta^{-\nu q}+e^{-\Gamma_{\varepsilon}^{\kappa}%
}\zeta^{-2\nu q}). \label{a11}%
\end{equation}
See the proof of Theorem 4.1 in \cite{[BF]} for detailed computations.

Moreover, in Proposition 4.4 in \cite{[BF]} one denotes $d_{t}=\det
\sigma_{F_{t}^{\varepsilon,\zeta}}$\ and proves that
\[
E(d_{t}^{-p})\leq C_{p,t}e^{c_{p}\Gamma_{\varepsilon}^{\gamma}}.
\]
Here $c_{p}$ is a constant which depends on $p$ and $C_{p,t}$ depends on $p$
but also on $t.$ The dependence in $t$ is not specified there, so we will
check here that
\begin{equation}
C_{p,t}=ct^{-4p\times\frac{2+\nu}{\nu}}. \label{a10}%
\end{equation}
We go in the proof of Proposition 4.4 in \cite{[BF]} and we find the
inequality
\[
E(d_{t}^{-p})\leq C_{p}e^{C_{p}\Gamma_{\varepsilon}^{\gamma}}\left(  \int
_{\xi\in R^{2}}\left\vert \xi\right\vert ^{8p-2}\exp(-ct\left\vert
\xi\right\vert ^{\nu/(2+\nu)})d\xi\right)  ^{1/2}.
\]
Using the change of variable $\overline{\xi}=t^{(2+\nu)/\nu}\xi$ one obtains
(\ref{a10}) so that
\[
\left\Vert d_{t}^{-1}\right\Vert _{p}\leq\frac{1}{t^{\frac{4(2+\nu)}{\nu}}%
}C_{p}e^{C_{p}\Gamma_{\varepsilon}^{\gamma}}.
\]
This, together with (\ref{a11}) and Schwarz's inequality gives (\ref{a13}).
$\square$


\bigskip

\addcontentsline{toc}{section}{References}

\begin{thebibliography}{99}                                                                                               %


\bibitem {ADVW}\textsc{R. Alexandre, L. Desvillettes,C. Villani,B. Wennberg}
(2000). Entropy dissipation and long-range interactions,
\emph{Arch.Rat.Mech.Anal.} \textbf{152}, 327-355.

\bibitem {BCpreprint}\textsc{V. Bally, L. Caramellino} (2017). Convergence and
regularity of probability laws by using an interpolation method. \emph{Ann..
Probab.} \textbf{45 }1110-1159 Preprint \texttt{arXiv:1409.3118}.

\bibitem {BCC}\textsc{V. Bally, L. Caramellino, R. Cont} (2016): Stochastic
integration by parts and functional It\^{o} calculus. \emph{Advanced Courses
in Mathematics - CRM Barcelona}\textit{,} Birkh\"{a}user.

\bibitem {[BCl1]}\textsc{V. Bally, E. Clement} (2011). Integration by parts
formulas and applications to equations with jumps. \emph{Probab. Theory
Related Fields} \textbf{151}, 613-657.

\bibitem {[BF]}\textsc{V. Bally, N. Fournier} (2011). Regularization
properties of the 2D homogeneous Boltzmann equation without cutoff.
\emph{Probab. Theory Related Fields} \textbf{151}, 659-704.

\bibitem {[BS]}\textsc{C. Bennett, R. Sharpley} (1988). Interpolation of
operators. Academic Press INC.

\bibitem {[BGJ]}\textsc{K. Bichteler, J.B. Gravereau J. Jacod }(1987).
Malliavin calculus for processes with jumps. \emph{Gordon and Breach science
publishers}, New York.

\bibitem {[Bi]}\textsc{J.M. Bismut} (1983). Calcul des variations stochastique
et processus de sauts. \emph{Z. Wahrsch. Verw. Gebiete} \textbf{63}, 147--235.

\bibitem {[BD]}\textsc{N. Bouleau, L. Denis}(2015). Dirichlet forms and
methods for Poisson point measures and L\'{e}vy processes. \emph{Probability
Theory and Stochastic Modelling}, \textbf{76}, Springer.

\bibitem {DF}\textsc{A. Debussche, N. Fournier} (2013). Existence of densities
for stable-like driven SDE's with Holder continuous coefficients. \emph{J.
Funct. Anal.} \textbf{264}, 1757-1778.

\bibitem {DR}\textsc{A. Debussche, M. Romito} (2014). Existence of densities
for the 3D Navier--Stokes equations driven by Gaussian noise. \emph{Probab.
Theory Related Fields} \textbf{158}, 575-596.

\bibitem {[De]}\textsc{S. De Marco} (2011). Smoothness and Asymptotic
Estimates of densities for SDEs with locally smooth coefficients and
Applications to square-root diffusions. \emph{Ann. Appl. Probab.} \textbf{21}, 1282-1321.

\bibitem {F3}\textsc{N. Fournier} (2000). Existence and regularity study for
$2D$ Bolzmann equation without cutoff by a probabilistic approach,. \emph{Ann.
Appl. Probab.} \textbf{10}, 434-462.

\bibitem {Ff}\textsc{N. Fournier} (2002). Jumping SDE's: absolute continuity
using monotonicity. \emph{Stochastic Process. Appl.} \textbf{98}, 317-330.

\bibitem {[F]}\textsc{N. Fournier} (2008). Smoothness of the law of some
one-dimensional jumping SDE's with non constant rate of jump. \emph{Electron.
J. Probab.} \textbf{13}, 135-156.

\bibitem {F1}\textsc{N. Fournier} (2015). Finiteness of entropy for the
homogeneous Boltzmann equation with measure initial condition. \emph{Ann.
Appl. Probab.} \textbf{25}, 860-897.

\bibitem {FM}\textsc{N. Fournier, C. Mouhot} (2009). On the well-posenss of
the spatially homogenous Bolzmann equation with a moderate angular singularity
\emph{Comm. Math. Phys.} \textbf{289}, 803-824.

\bibitem {[FP]}\textsc{N. Fournier, J. Printems} (2010). Absolute continuity
of some one-dimensional processes. \emph{Bernoulli} \textbf{16}, 343-360.

\bibitem {GPV}\textsc{I.M. Gamba, V. Panferov, C. Villani} (2009). Upper
Maxwellian bounds for the spacialy homogenuous Bolzmann equation \emph{Arch.
Ration. Mech. Anal. } \textbf{194}, 253-282.

\bibitem {GM}\textsc{C. Graham,S Meleard} (1999). Existance and regularity of
a solution of a Kac equation without cutoff using the stochastic calculus of
variations. \emph{Comm. Math. Phys..} \textbf{205,} \textbf{no 3}, 551-569.

\bibitem {GMN}\textsc{H. Guerin,S Meleard,E. Nualart} (2006). Exponential
estimates for spatially homogeneous Landau equations via the Malliavin
calculus. \emph{Journal of Func. Analysis, }\textbf{no. 2}, 649--677.

\bibitem {[I]}\textsc{Y. Ishikawa }(2013). \emph{Stochastic Calculus of
variation for Jump Processes}. De Gruyter Studies in Math. \textbf{54}.

\bibitem {K}\textsc{V.N. Kolokoltsov} (2006). On the regularity of the
solutions of trhe space homogenuous Bolzmann equation with polynomial growing
collision kernel. \emph{Advanced Studies in Contemporany Math. , }\textbf{no.
21}, 9--38.

\bibitem {[L]}\textsc{R. Leandre} (1985). R\'{e}gularut\'{e} des processus de
sauts g\'{e}g\'{e}n\'{e}r\'{e}s. \emph{Ann. Inst. H. Poincar\'{e} Probab.
Statist.} \textbf{21}, 125--146.

\bibitem {PSL}\textsc{E. Pardoux,R Sentis,B. Lapeyre.} (2003) Introduction to
Monte-Carlo Methods for Transport and Diffusion Equations \emph{Oxford texts
in applied and engineering mathematics, }\textbf{6.}

\bibitem {T}\textsc{H. Tanaka} (1978). Probabilistic treatment of the
Bolzmaznn equation of Maxwellian molecules \emph{Z,Wahrsch. und Verw. Gebiete
}\textbf{46}, 67--105.
\end{thebibliography}
\end{document}